\documentclass[12pt]{article}

\usepackage{graphicx}
\usepackage{arydshln}
\usepackage{amssymb}

\begin{document}
\centerline{A modified Conway--Maxwell--Poisson type binomial distribution and its applications}
\vskip 3mm

\vskip 5mm
\noindent T. Imoto$^a$, C. M. Ng$^b$, S. H. Ong$^b$ and S. Chakraborty$^c$

\noindent $^a$ School of Management and Information, University of Shizuoka, 51-1 Yada, Suruga-ku, Shuzuoka 422-8526, Japan

\noindent $^b$ Institute of Mathematical Sciences, Faculty of Science, University of Malaya, 50603 Kuala Lumpur, Malaysia

\noindent $^c$ Department of Statistics, Dibrugarh University, Dibrugarh-786004, Assam, India

\noindent Corresponding author E-mail: imoto0923@gmail.com

\vskip 3mm

\noindent {\bf Abstract}
This paper proposes a generalized binomial distribution with four parameters, 
which is derived from the finite capacity queueing system with state-dependent service and arrival rates. 
This distribution is also generated from the conditional Conway--Maxwell--Poisson distribution 
given a sum of two Conway--Maxwell--Poisson variables.
In this paper, we consider the properties about the probability mass function, 
index of dispersion, skewness and kurtosis and give applications of the proposed distribution from its geneses.  
The estimation method and simulation study are also considered.
\vskip 4mm

\noindent {\bf 1. Introduction}

The binomial distribution is classically utilized as a model for analyzing count data with finite support
since it has a simple genesis arising from Bernoulli trials and belongs to the exponential family that can model under-dispersion, 
in which the variance is smaller than the mean.  
However, its reliance on the simple genesis limits its flexibility in many applications.
For examples, the binomial distribution cannot model over-dispersion, 
lepto-kurtosis when the mode is around $n/2$ and platy-kurtosis when the mode is around $0$ or $n$, 
where $n$ is an integer parameter of the binomial distribution.
For overcoming the limitations, many researchers have developed various generalized binomial distributions 
through considering more practical assumption and adding more parameters.

Consul (1974) and Consul and Mittal (1975) derived generalized binomial distributions, referred to the quasi-binomial distributions,
by considering the urn problems where a person decides his strategy before making draws from the urn.
Detailed studies about the distributions have been summarized in Consul and Famoye (2006).
There are also many generalizations from the sum of dependent Bernoulli random variables.
For example, Altham (1978) used the multiplicative and additive binary variables whereas
Chang and Zelterman (2002) used the binary variables whose conditional probabilities of success given the previous status 
depend only on the number of successes of previous variables. 
Upton and Lampitt (1981) considered a convolution of the binomial and Poisson variables 
to describe the changes in the counts of bird territories in successive years and
Ong (1988) independently considered the same distribution, 
which was referred as the Charlier series distribution, and gave detailed studies about its properties and characteristics.
The generalized Charlier series distribution by Kitano et al. (2005) was derived from two-step recursive formula and applied to the collective risk theory.

Conway and Maxwell (1962) considered the queueing models with state-dependent service or arrival rate 
and derived a generalized Poisson distribution with probability mass function (pmf)
\begin{eqnarray*}
\label{compmf}
{\rm P}(X=x)=\frac{\lambda^x}{x!^r}\frac{1}{Z(\lambda,r)},\quad \mbox{where}\quad Z(\lambda,r)=\sum_{k=0}^{\infty}\frac{\lambda^k}{(k!)^r}
\end{eqnarray*}
for $r>0$ and $\lambda>0$, which is known as the Conway--Maxwell--Poisson (CMP) distribution. 
In this paper, we write $X \sim$ CMP$(r ,\lambda)$ to indicate that the random variable $X$ has the pmf (\ref{compmf}).
The CMP distribution has been revived by Shmueli et al. (2005) as a flexible distribution 
that can adapt to over-dispersion for $r<1$ or under-dispersion for $r>1$. 
Considering the conditional CMP distribution given a sum of two CMP variables 
with the same dispersion parameter and different central parameters,
Shmueli et al. (2005) defined the Conway--Maxwell--Poisson type binomial (CMPB) distribution 
and Borges et al. (2014) gave detailed studies about this distribution.

In this paper, we consider a generalized binomial distribution with four parameters, which is derived from 
a finite capacity queueing system with state-dependent arrival and service rates.
The distribution is also derived from the conditional CMP distribution given a sum of two CMP variables 
with the different dispersion and central parameters.
From these geneses, this distribution may have biological and marketing applications 
such as modeling the number of males in a given number of group, items in a storage and individuals in a colony. 
The proposed distribution is appealing from a theoretical point of view since it belongs to the exponential family. 
Various statistical and probabilistic properties can be derived such as moments and estimations. 
Furthermore, the proposed distribution includes many interesting distributions; 
a degenerate, uniform, CMPB, truncated CMP and new Conway--Maxwell--Poisson type generalized binomial distributions.
The advantages of this distribution is its flexibility of the dispersion, skewness, kurtosis and modality. 
This distribution can become a under- or over-dispersed distribution, positively or negatively skewed distribution, 
lepto- or platy-kurtic distribution as well as unimodal or bimodal distribution. 
These versatility and flexibility give well performance for various datasets.

In this paper, Section 2 defines a generalized binomial distribution with four parameters 
and considers its genesis arising from the finite capacity queueing process.
In the same section, we gives the relationship between the derived distribution and the CMP distribution.
Examples of applications are given in this section.
Section 3 considers the distributional properties and introduce some special distributions.  
Section 4 discusses the computational aspect of the proposed distribution while
Section 5 considers the maximum likelihood (ML) estimation and the simulation study to show the performance of the estimation.
Using the ML estimation, we give illustrative examples by fitting the proposed distribution to real datasets in Section 6.
Finally, concluding remarks are given in Section 7.

\vskip 3mm

\noindent {\bf 2. Modified Conway--Maxwell--Poisson type binomial distribution}

\noindent {\bf 2.1. Definition}

Let us consider the distribution with pmf 
\begin{eqnarray}
\label{MCMPBpmf}
{\rm P}(X=x)=\frac{\theta^x}{{x!}^{\alpha}{(n-x)!}^{\beta}C_n(\alpha,\beta,\theta)} ,\ x=0,1\ldots,n,
\end{eqnarray}
for $\alpha, \beta \in (-\infty, \infty)$, $\theta>0$ and a positive integer $n$, where 
\begin{eqnarray*}
\label{MCMPBconstant}
C_n(\alpha,\beta,\theta)=\sum_{k=0}^n\frac{\theta^k}{{k!}^{\alpha}{(n-k)!}^{\beta}}.
\end{eqnarray*}
This distribution reduces to the CMPB distribution by Shmueli et al. (2005) and Borges et al. (2014) when $\alpha=\beta$.
Hence we shall refer to distribution (\ref{MCMPBpmf}) as the modified Conway--Maxwell--Poisson type binomial (MCMPB) distribution 
and denote $X \sim$ MCMPB$_n(\alpha,\beta,\theta)$.
The following subsections explain the genesis of the MCMPB distribution, and the relationship between the MCMPB and CMP distributions.

\vskip 3mm

\noindent {\bf 2.2.  Finite capacity queueing process}

Consider a single-queue-single-server system,
where the service times are exponentially distributed with mean $\mu x^{\alpha}$ and 
the customers arrive according to a Poisson process with mean $\lambda(n-x)^{\beta}$ 
for $x\leq n$ and $0$ otherwise when the size of the queue is $x$.
Customers are served on first-come-first-served basis.
In this system, the capacity of this queue is finite ($=n$), and 
the service and arrival rates are increasing (decreasing) as the size of queue is increasing 
when $\alpha>(<)0$ and $\beta<(>)0$. 

The probability $P(x,t)$ when the size of the queue is $x$ at time $t$ satisfies the following difference equation for small $h>0$
\begin{eqnarray*}
P(x,t+h)=\mu (x+1)^{\alpha}h P(x+1,t)+\{1-\mu x^{\alpha}h-\lambda(n-x)^{\beta}h \}P(x,t)\\
+\lambda(n-x+1)^{\beta}h P(x-1,t)
\end{eqnarray*}
for $x=0,1,\ldots,n$ with the condition $P(-1,t)=P(n+1,t)=0$.
Letting $h\rightarrow0$ in the above equation, we get the following difference-differential equation
\begin{eqnarray*}
\frac{\partial P(x,t+h)}{\partial t}=\mu(x+1)^{\alpha}P(x+1,t)-\{\mu x^{\alpha}+\lambda(n-x)^{\beta}\}P(x,t)\\
+\lambda(n-x+1)^{\beta} P(x-1,t).
\end{eqnarray*}
Assuming a stationary state, or $\partial P(x, t)/\partial t=0$, and then putting $P(x,t)=P(x)$ and $\theta=\lambda/\mu$, 
we finally have the difference equation
\begin{eqnarray*}
P(x+1)=\frac{x^{\alpha}+\theta(n-x)^{\beta}}{(x+1)^{\alpha}}P(x)-\frac{\theta(n-x+1)^{\beta}}{(x+1)^{\alpha}}P(x-1),
\end{eqnarray*}
for $x=0,1,\ldots,n$ with the condition  $P(-1)=P(n+1)=0$.
The solution of this equation is proved by induction to be 
$$
P(x)=\frac{n!^{\beta}P(0) \theta^x}{x!^{\alpha}(n-x)!^{\beta}} ,\quad x=0,1\ldots,n.
$$
We can see that $\{n!^{\beta}P(0)\}^{-1}=C_n(\alpha,\beta,\theta)$ and this distribution is equal to MCMPB$_n(\alpha,\beta,\theta)$.

From this genesis, we can apply the MCMPB distribution to the size of water unit in dam or commodities in storage 
since the capacity of dam or storage is usually finite.  
Moreover, treating the arrival and service rates as the birth and death rates respectively, 
we can consider this process as the birth and death process.
This interpretation leads to the applications in biology such as modeling the size of individuals in some colony.\\

\noindent{\it Example }: Bacterial clumps in a milk film (Bliss 1953).

In this dataset, a microscope slide was split into 400 regions of equal area and 
the number of bacterial clumps on each was counted.
When we fit the MCMPB$_n(\alpha,\beta, {\rm e}^{\psi})$ to this dataset 
by the method introduced in Section 5,
the profile maximum likelihood estimates of $n$ is $19$ and maximum likelihood estimates and 
$95\%$ confidence intervals of $\alpha$, $\beta$, $\psi$ are 
$0.73$ and $(0.55, 0.92)$, $-1.00$ and $(-1.34, -0.66)$, $3.35$ and $(2.22, 4.47)$, 
respectively (The $p$-value of $\chi^2$ test is $0.35$).
From this result, we see that the capacity of each region is $19$ and 
both birth and death rates of bacterial clumps increase as the number of bacterial clumps in a field increases
because the confidence intervals confirm $\alpha>0$ and $\beta<0$.
\begin{table}[htbp]
\caption{The number of bacterial clumps per field in a milk film.}
\label{bacterial}
\begin{center}
{\scriptsize \begin{tabular}{l}
\begin{tabular}{c|cccccccccc}
Count    & 0     & 1     & 2     & 3     & 4     & 5     & 6     & 7    & 8    & 9       \\
\hline
Observed & 56    & 104   & 80    & 62    & 42    & 27    & 9     & 9    & 5    & 3        \\
Expected & 60.65 & 91.01 & 86.79 & 65.14 & 42.07 & 24.62 & 13.51 & 7.13 & 3.70 & 1.92  \\
\end{tabular}
\\
\\
\begin{tabular}{c|ccccccccccc}
Count    & 10 & 11 & 12 & 13 & 14 & 15 & 16 & 17 & 18 &19 & Total  \\
\hline
Observed & 2 & 0    & 0    & 0    & 0    & 0    & 0    & 0    & 0    & 1    & 400 \\
Expected & 1.01 & 0.56 & 0.32 & 0.20 & 0.14 & 0.11 & 0.10 & 0.12 & 0.21 & 0.69 & 400.00 \\
\end{tabular}
\end{tabular}}
\end{center}
\end{table}

\vskip 3mm

\noindent {\bf 2.3. Conditional CMP distribution}

For the MCMPB distribution, we have the following theorem which extends the relationship between the Poisson and binomial distributions. \\

\noindent {\it Theorem 1.}\\
If the random variables $X_1$ and $X_2$ are independently distributed as the CMP distributions,
then the conditional distribution of $X_1$ given $X_1+X_2=n$ is distributed as the MCMPB distribution.
Conversely, if the conditional distribution of $X_1$ given $X_1+X_2=y$ is distributed as the MCMPB$_y(\alpha, \beta, \theta)$ for $y=0,1,\ldots,$
and $X_1$ is independent of $X_2$, then $X_1$ and $X_2$ are distributed as the CMP distributions.\\

\noindent {\it Proof.}
Assume that $X_1\sim$ CMP$(\alpha, \lambda_1)$, $X_2\sim$ CMP$(\beta,\lambda_2)$ 
and these random variables are independent.
The conditional distribution of $X_1$ given $X_1+X_2=n$ is easily obtained as
\begin{eqnarray*}
{\rm P}(X_1=x|X_1+X_2=n)
&=&\frac{ {\rm P}(X_1=x) {\rm P}(X_2=n-x) }{\sum_{k=0}^n   {\rm P}(X_1=k) {\rm P}(X_2=n-k) }\\
&=&\frac{(\lambda_1/\lambda_2)^x}{{x!}^{\alpha}{(n-x)!}^{\beta}C_n(\alpha,\beta,\lambda_1/\lambda_2)}.
\end{eqnarray*}
This is equal to the pmf of  MCMPB$_n(\alpha,\beta,\lambda_1/\lambda_2)$.

Next we prove the second statement.
Assume that the pmf of $X_1|X_1+X_2=y$ is that of MCMPB$_y(\alpha, \beta, \theta)$ for $y=0,1,\ldots$.
Then we see
\begin{eqnarray}
\label{conditionalproof}
\frac{ {\rm P}(X_1=x) {\rm P}(X_2=y-x) }{ {\rm P}(X_1=x-1) {\rm P}(X_2=y-x+1) } 
&=& \frac{ {\rm P}(X_1=x|X_1+X_2=y) }{ {\rm P}(X_1=x-1|X_1+X_2=y) } \nonumber\\
&=& \frac{ (y-x+1)^{\beta} \theta }{ x^{\alpha} }.
\end{eqnarray}
Put $y=x$ in (\ref{conditionalproof}), we get
$$
{\rm P}(X_1=y)  = \frac{ a \theta }{ y^{\alpha}}{\rm P}(X_1=y-1)=\cdots=\frac{(a\theta)^y}{y!^{\alpha}}{\rm P}(X_1=0),
$$
where $a={\rm P}(X_2=1)/{\rm P}(X_2=0)$. From this equation, we can see that $X_1 \sim$ CMP$(\alpha, a\theta)$.
Similarly, letting $x=1$ in (\ref{conditionalproof}), we obtain
$$
{\rm P}(X_2=y)  = \frac{{\rm P}(X_2=y-1) } {b \theta  y^{\beta}} = \cdots = \frac{{\rm P}(X_2=0) } {(b \theta)^y  y!^{\beta}}
$$
where $b={\rm P}(X_1=0)/{\rm P}(X_1=1)$. 
From the equation (\ref{conditionalproof}), we see that $(ab)^{-1}=\theta$ and thus, $X_2 \sim$ CMP$(\beta, a)$. 
\begin{flushright}$\Box$\end{flushright}

For $\alpha=\beta$, the first statement in Theorem 1
is the genesis of the CMPB distribution by Shmueli et al. (2005) and, 
in this sense, the proposed distribution modifies the CMPB distribution. 

As an example of applications of Theorem 1, 
we can consider the problem of modeling the number of males in some group.
When we analyze the trend of the number of males, 
the sampling for $X_1|X_1+X_2=n$ might be easier than the direct sampling for $X_1$,  
where the random variables $X_1$ and $X_2$ represent the numbers of males and females respectively.
By fitting the MCMPB distribution to the the sampling data for $X_1|X_1+X_2=n$, 
we can see the trend of males through the CMP distribution. \\

\noindent {\it Example}: Male children in family (Sokal and Rohlf 1994) (cf. Lindsey 1995).

In this dataset, the number of male children in 6115 families with 12 children 
each in the nineteenth century Saxony was counted.
When we fit the MCMPB$_n(\alpha,\beta, {\rm e}^{\psi})$ to this dataset by the method introduced in Section 5,
the maximum likelihood estimates and $95\%$ confidence intervals of $\alpha$, $\beta$, $\psi$ are 
$0.93$ and $(0.74, 1.12)$, $0.76$ and $(0.59, 0.94)$, $0.37$ and $(-0.28, 1.04)$, respectively 
(The $p$-value of $\chi^2$ test is $0.13$).
From these results with Theorem 1,
the number of males in a family in Saxony can be said to be distributed as the Poisson distribution (since $\alpha \approx 1.0$)
while that of females was distributed as the over-dispersed CMP distribution (since $\beta < 1.0$).
\begin{table}[htbp]
\caption{The number of males in 6115 families with 12 children in Saxony.}
\label{children}
\begin{center}
{\scriptsize
\begin{tabular}[c]{l}
\begin{tabular}[c]{c|ccccccccccccccccccccccccccccccc}
Count    & 0    & 1     & 2      & 3      & 4      & 5       & 6     \\
\hline
Observed & 3    & 24    & 104    & 286    & 670    & 1033    & 1343     \\
Expected & 2.22 & 21.49 & 102.00 & 308.64 & 659.30 & 1045.91 & 1264.63 \\
\end{tabular}\\
\\
\begin{tabular}[c]{c|ccccccccccccccccccccccccccccccc}
Count    & 7       & 8      & 9      & 10     & 11    & 12   & Total  \\
\hline
Observed & 1112    & 829    & 478    & 181    & 45    & 7    & 6115 \\
Expected & 1177.77 & 842.95 & 456.07 & 179.65 & 47.54 & 6.84 & 6115.00 \\
\end{tabular}
\end{tabular}}
\end{center}
\end{table}

\vskip 3mm

\noindent {\bf 3. Properties of MCMPB distribution}

\noindent {\bf 3.1. Exponential family}

From the expression of pmf (\ref{MCMPBpmf}), the MCMPB$_n(\alpha, \beta, \theta)$ 
belongs to the power series distributions with power parameter $\theta$
and thus, has the recursive formula about the moments as
\begin{eqnarray*}
\label{MCMPBmoments}
\left\{\begin{array}{l}
\displaystyle  \mu'_{1}=\theta \frac{\partial \log C_n(\alpha,\beta,\theta)}{\partial \theta},\\
\displaystyle \mu'_{k+1}=\theta \frac{\partial \mu'_k}{\partial \theta}+ \mu'_{1}\mu'_k,\\
\displaystyle \mu_{k+1}=\theta \frac{\partial \mu_k}{\partial \theta}+k\mu_2\mu_{k-1},
\end{array}\right.
\end{eqnarray*}
where $\mu'_k={\rm E}[X^{k}]$ and $\mu_k={\rm E}[(X-{\rm E}[X])^{k}]$ 
with $X \sim$ MCMPB$_n(\alpha, \beta, \theta)$.

Moreover, by replacing $\theta= {\rm e}^{\phi}$ in (\ref{MCMPBpmf}), 
we can rewrite the pmf of the MCMPB distribution as
\begin{eqnarray}
\label{MCMPBpmf2}
{\rm P}(X=x)=\exp\{x \phi-\alpha \log x!-\beta \log(n-x)!-\log C^*_n(\alpha,\beta,\phi)\},
\end{eqnarray}
where $C^*_n(\alpha,\beta,\phi)=C_n(\alpha,\beta,{\rm e}^{\phi})$ for a positive integer $n$,
and $\alpha, \beta, \phi \in (-\infty, \infty)$.
From the expression (\ref{MCMPBpmf2}), we see that the MCMPB distribution is a member of 
the exponential family with natural parameters $(\phi, \alpha, \beta)$
when $n$ is known.
From this expression, the first and second moments about $X$, $\log X!$, $\log (n-X)!)$ are obtained as
\begin{eqnarray}
\label{covariance}{\footnotesize
\left\{\begin{array}{l}
\displaystyle {\rm E}[X]=\frac{\partial \log C^*_n}{\partial \phi},\ 
\displaystyle {\rm Var}[X]=\frac{\partial^2 \log C^*_n}{\partial \phi^2},\
\displaystyle {\rm E}[\log X!]=-\frac{\partial \log C^*_n}{\partial \alpha},\\ 
\displaystyle {\rm Var}[\log X!]=\frac{\partial^2 \log C^*_n}{\partial \alpha^2},\
\displaystyle {\rm E}[\log (n-X)!]=-\frac{\partial \log C^*_n}{\partial \beta},\\\
\displaystyle {\rm Var}[\log (n-X)!]=\frac{\partial^2 \log C^*_n}{\partial \beta^2},\
\displaystyle {\rm Cov}[X, \log X!]=-\frac{\partial^2 \log C^*_n}{\partial \phi \partial \alpha},\\ 
\displaystyle {\rm Cov}[\log X!,\log (n-X)!]=\frac{\partial^2 \log C^*_n}{\partial \alpha \partial \beta},\
\displaystyle {\rm Cov}[\log (n-X)!,X]=-\frac{\partial^2 \log C^*_n}{\partial \beta \partial \phi}.\\
\end{array}\right.}
\end{eqnarray}

\vskip 3mm

\noindent {\bf 3.2. Exponential combination}

MCMPB$_n(\alpha,\beta,\theta)$ is an exponential combination of binomial distribution with parameters $(n,p)$
and truncated version of CMP$(\lambda, r)$, with $\theta=\lambda\{p/\lambda(1-p)\}^{\beta}$ and $\alpha=\beta(1-r)+r$.
From this fact, it is clear that, for MCMPB$_n(\alpha,\beta,\theta)$, $\beta$ close to zero indicates departure from the binomial distribution 
towards truncated the CMP distribution, while $\beta$ close to one indicates the reverse.
In case the value $\beta$ close to $1/2$ will indicate that both the distribution fit the data equally well.
Thus MCMPB$_n(\alpha,\beta,\theta)$ can also be regarded as a natural extension of the binomial and truncated CMP distributions.

\vskip 3mm

\noindent {\bf 3.3. Sum of dependent Bernoulli variables}

Let $(X_1,X_2,\ldots,X_n)$ be the Bernoulli variables having joint pmf 
$$
{\rm P}(X_1=x_1,\ldots,X_n=x_n)=\frac{p^{\bf x}(1-p)^{n-{\bf x}}}{{\bf x}!^{\alpha-1}(n-{\bf x})!^{\beta-1}} 
\left/\sum_{k_i's=0, 1}  \frac{p^{\bf k}(1-p)^{n-{\bf k}}}{{\bf k}!^{\alpha-1}(n-{\bf k})!^{\beta-1}}   \right.
$$
for $0<p<1$, where ${\bf x}=x_1+\cdots+x_n$ and ${\bf k}=k_1+\cdots+k_n$.
Then the distribution of $X=\sum_{i=0}^n X_i$ is
\begin{eqnarray*}
{\rm P}(X=x)&=&\sum_{x_1+\cdots+x_n=x} {\rm P}(X_1=x_1,\ldots, X_n=x_n)\\
                      &=& {n \choose x} \frac{\theta^{x}}{x!^{\alpha-1}(n-x)!^{\beta-1}} 
                             \left/\sum_{k=0}^n  {n \choose k} \frac{\theta^{k}}{k!^{\alpha-1}(n-k)!^{\beta-1}},   \right. \\
                      &=& \frac{\theta^{x}}{x!^{\alpha}(n-x)!^{\beta}C_n(\alpha,\beta,\theta)}, 
\end{eqnarray*}
where $\theta=p/(1-p)$. This is the pmf of MCMPB$_n(\alpha,\beta,\theta)$.

Here, by assuming that, for $X_1, \ldots,_n$, all correlations higher than order two are zero (Bahadur 1961),
we can derive the expectation and correlation between any pairs of $(x_i, X_j)$ ($i \not = j = 1$) as
\begin{eqnarray*}
{\rm E}[X_i]=\frac{2^{\alpha-1}p(1-p)+p^2}{2^{\alpha -\beta}(1-p)^2+2^{\alpha}p(1-p)+p^2}
=\frac{2^{\alpha-1}\theta+\theta^2}{2^{\alpha -\beta}+2^{\alpha}\theta+\theta^2},
\end{eqnarray*}
\begin{eqnarray*}
{\rm Cor}[X_i, X_i]&=&\frac{p(1-p)(2^{\alpha-\beta}-4^{\alpha-1})}{\{2^{\alpha-1}(1-p)+p\}\{2^{\alpha-\beta}(1-p)+2^{\alpha-1}p\}}\\
&=&\frac{\theta(2^{\alpha-\beta}-4^{\alpha-1})}{(1+\theta)\{2^{\alpha-1}+\theta\}\{2^{\alpha-\beta}+2^{\alpha-1}\theta\}}.
\end{eqnarray*}
These are the extension of the result by Borges et al. (2014).

\vskip 3mm

\noindent {\bf 3.4. Behavior of pmf}

When $X \sim$ MCMPB$_n(\alpha, \beta, \theta)$, we see that
$$
\frac{{\rm P}(X=x+1){\rm P}(X=x-1)}{{\rm P}(X=x)^2}=\left(\frac{x}{x+1}\right)^{\alpha} \left(\frac{n-x}{n-x+1}\right)^{\beta}
$$
and the pmf has log-concavity for $\alpha, \beta>0$ and the log-convexity for $\alpha, \beta<0$. 
From this fact, we see that the MCMPB distribution with $\alpha, \beta>0$ is always unimodal
and has strong unimodality, the property that its convolution with any unimodal distribution is unimodal.
Note that the mode is not at $x=0$ or $x=n$ in this case because ${\rm P}(X=1)/{\rm P}(X=0)=n^{\beta}>1$ 
and ${\rm P}(X=n)/{\rm P}(X=n-1)=n^{-\alpha}<1$.
Also, we see that the MCMPB distribution with $\alpha, \beta<0$ is a unimodal distribution 
at the mode $x=0$ when $\theta<n^{\alpha}$ or at the mode $x=n$ when $\theta>n^{-\beta}$
and a bimodal distribution at the modes $x=0$ and $x=n$ when $n^{\alpha}<\theta<n^{-\beta}$.
For either $\alpha<0$ or $\beta<0$, the MCMPB distribution might become a bimodal distribution whose one mode is at $x=0$ or $x=n$.

\begin{figure}[h]
\begin{center}
\begin{tabular}{cc}
\includegraphics[width=0.45\hsize,clip]{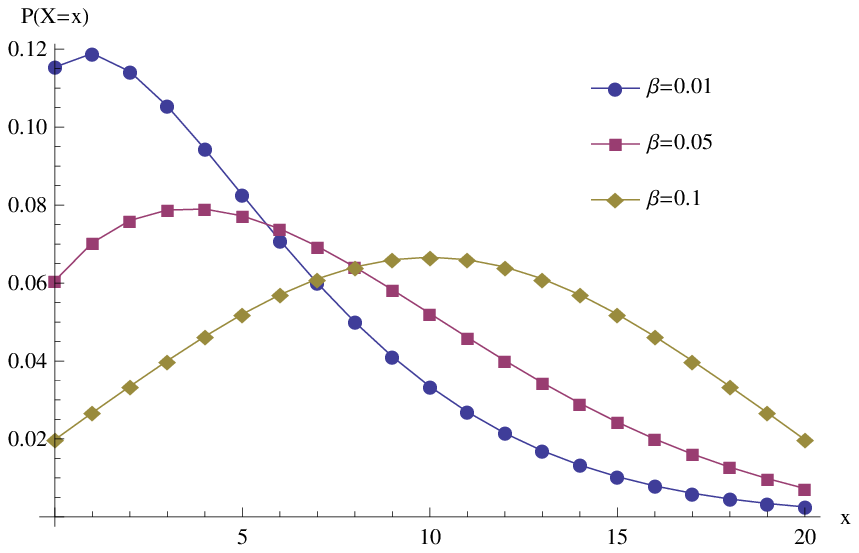}
& 
\includegraphics[width=0.45\hsize,clip]{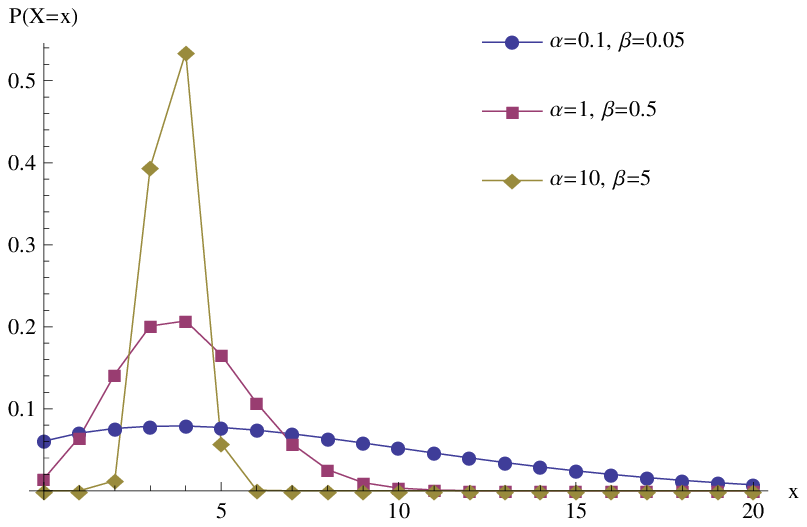}\\
(a) MCMPB$_{20} (0.1, \beta, 1)$ & (b) MCMPB$_{20}(\alpha, \beta, 1)$\\
\includegraphics[width=0.45\hsize,clip]{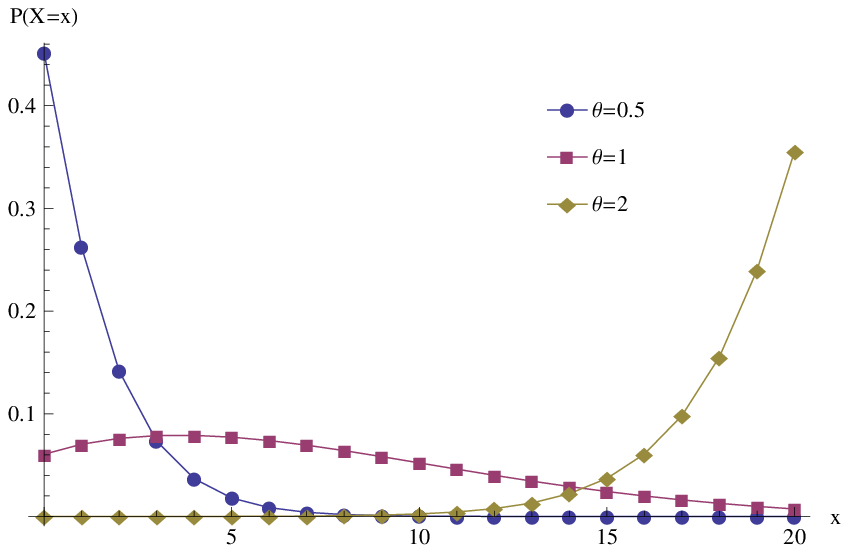} 
& 
\includegraphics[width=0.45\hsize,clip]{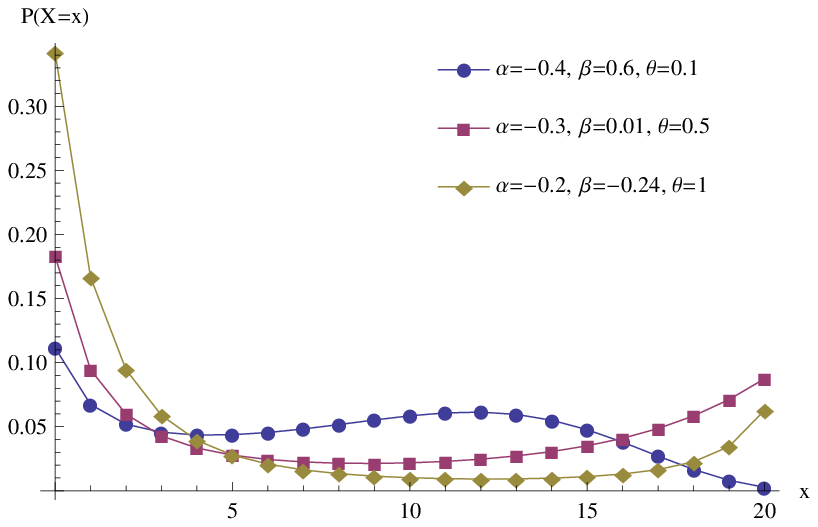}\\
(c) MCMPB$_{20}(0.1, 0.05, \theta)$ &(d) MCMPB$_{20}(\alpha, \beta, \theta)$ 
\end{tabular}
\caption{The plot of the pmf of the MCMPB distribution.}
\label{pmf}
\end{center}
\end{figure}

The graphical illustrations for the pmf of the MCMPB distribution are given in Figure \ref{pmf}.
Figure \ref{pmf} (a) reveals that the mode shifts from the left to the right as $\beta>0$ becomes large
while Figure \ref{pmf} (b) shows that, for $\theta=1$, the MCMPB distribution becomes a platy-kurtic distribution when $\alpha,\beta>0$ are small 
and a lepto-kurtic distribution when $\alpha,\beta>0$ are large.
However, Figure \ref{pmf} (c) shows that, for $\theta\not=1$, the distribution becomes a lepto-kurtic distribution even when  $\alpha,\beta>0$ are small.
This is because MCMPB$_n(\alpha,\beta,1)$ approaches to a uniform distribution when $\alpha,\beta \rightarrow 0$
whereas MCMPB$_n(\alpha,\beta,\theta)$ with $0<\theta<1$ approaches to a truncated geometric distribution when $\alpha,\beta \rightarrow 0$.
Interestingly, the MCMPB distribution can become a platy-kurtic distribution with any mode 
whereas the CMPB distribution becomes a platy-kurtic distribution only when the mode is at $x=n/2$.
Figure \ref{pmf} (d) shows that the MCMPB distribution might become a bimodal distribution, whose one mode is at $x=0$ or $x=n$ when $\alpha<0$ or $\beta<0$.

It is seen that $X\sim$ MCMPB$_n(\alpha,\beta,\theta)$ $\Leftrightarrow$ $n-X\sim$ MCMPB$_n(\beta,\alpha,1/\theta)$
and this leads to the figures which are horizontal inversions of Figure \ref{pmf}.
This fact also leads to the formula
$$
{\rm P}(X \geq x|n,\alpha,\beta,\theta)={\rm P}(X\leq n-x|n,\beta,\alpha,1/\theta),
$$
which will be useful for calculating the distribution function for large $n$.

\vskip 3mm

\noindent {\bf 3.5. Special cases}

The MCMPB distribution includes some interesting distributions as follows.

\begin{itemize}

\item 
The degenerate distribution at $x=a$\\
This is proven as follow:
Let $f(x)=x!^{-\alpha}(n-x)!^{-\beta}$ for $x=0,1,\ldots,n$.
MCMPB$_n(\alpha,\beta,1)$ is unimodal when $\alpha, \beta>0$.
If the mode is at $x=a(\not=0,n)$, we have $f(a-1) < f(a) > f(a+1)$, 
or $\log(a)/\log(n-a+1) < \beta/\alpha < \log(a+1)/\log(n-a)$. 
Put $d=\beta/\alpha$.
This leads to $f(k_1)/f(a)<f(a-1)/f(a)=\{a/(n-a+1)^d\}^{\alpha} \rightarrow 0$ 
and $f(k_2)/f(a)<f(a+1)/f(a)=\{(n-a)^d/(a+1)\}^{\alpha} \rightarrow 0$ 
for $k_1=0,1,\ldots,a-2$ and $k_2=a+2,a+3,\ldots$ as $\alpha, \beta \rightarrow \infty$.
Therefore, we have the result $f(a)/C_n(\alpha,\beta,1) \rightarrow 1$ as $\alpha, \beta \rightarrow \infty$.
From the similar argument for $\alpha, \beta<0$ and $n^{\alpha}>1$ or $n^{-\beta}<1$,
We can see that MCMPB$_n(\alpha,\beta,1)$ reduces to the degenerate distribution at $x=0$ or $x=n$.  \\

\item The uniform distribution at $x=a,a+1$\\ 
This is proven by considering the case when the mode of MCMPB$_n(\alpha,\beta,1)$ is at $x=a,a+1$ in the above proof. \\

\item The discrete uniform distribution\\ 
This is equal to MCMPB$_n(0,0,1)$.

\end{itemize}

The MCMPB distribution includes the degenerate and uniform distributions, 
indicating that the MCMPB distribution can become a lepto- or platy-kurtic distribution. 

\begin{itemize}
\item 
The truncated CMP distribution 
$$
{\rm P}(X=x)=\frac{\theta^x}{{x!}^{\alpha}C_n(\alpha,0,\theta)},\ x=0,1\ldots,n.
$$
This is equal to MCMPB$_n(\alpha,0,\theta)$ with $\alpha>0$ and $0<\theta<1$.
The normalizing constant of the truncated CMP distribution includes only the summation of finite series.
Therefore, it is easier to apply this distribution than the CMP distribution.\\

\item The CMP distribution\\
This occurs when $n \rightarrow \infty$ and $n^{\beta}\theta=\lambda$ (finite) for MCMPB$_n(\alpha,\beta,\theta)$.
This is proven as follow: 
\begin{eqnarray*}
&  &\lim_{n \rightarrow \infty} \frac{\theta^x}{x!^{\alpha}(n-x)!^{\beta}} \left/ \sum_{k=0}^n \frac{\theta^k}{k!^{\alpha}(n-k)!^{\beta}}  \right.\\
&=&\lim_{n \rightarrow \infty} \frac{n!^{\beta}\theta^x}{x!^{\alpha}(n-x)!^{\beta}} \left/ \sum_{k=0}^n \frac{n!^{\beta}\theta^k}{k!^{\alpha}(n-k)!^{\beta}}  \right.\\
&=&\lim_{n \rightarrow \infty} \frac{\{n\cdots(n-x+1)\}^{\beta}\theta^x}{x!^{\alpha}} \left/ \sum_{k=0}^n \frac{\{n\cdots(n-k+1)\}^{\beta}\theta^k}{k!^{\alpha}}  \right.\\
&=&\frac{\lambda^x}{x!^{\alpha}}\left/ \sum_{k=0}^{\infty} \frac{\lambda^k}{k!^{\alpha}} \right.
\end{eqnarray*}\\

\item 
A dispersion distribution
$$
{\rm P}(X=x)=\frac{\theta^x}{x!{(n-x)!}^{\beta}C_n(1,\beta,\theta)},\ x=0,1\ldots,n.
$$
This is equal to MCMPB$_n(1,\beta,\theta)$.
This distribution corresponds to the weighted version of the Poisson distribution (Rao 1965) 
with the weight function $w(x)=\Gamma(n-x+1)^{-\beta}$ for $0 \leq x \leq n$ and $0$ otherwise.
Corollary 4 in Castillo and P\'{e}rez-Casany (2005) confirms that the weighted Poisson distribution with the weight function $w(x)=\exp\{rt(x)\}$, 
where $t(\cdot)$ is a convex function, is over-dispersed for $r>0$ and under-dispersed for $r<0$. 
Since $\log \Gamma(n-x+1)$ is a convex function, this distribution is over-dispersed for $\beta<0$ and under-dispersed for $\beta>0$.\\

\item
A skew distribution 
$$
{\rm P}(X=x)=\frac{1}{{x!}^{\alpha}{(n-x)!}^{\beta}C_n(\alpha,\beta,1)},\ x=0,1\ldots,n.
$$
This is equal to MCMPB$(n,\alpha,\beta,1)$. 
As we can see from Figure \ref{pmf} with $\theta=1$, the parameters $\alpha$ and $\beta$ controls the mean and variance of the distribution.
And, the indicator $\alpha-\beta$ measure the degree of skewness.

\end{itemize}

\vskip 3mm

\noindent {\bf 3.6. Indices of dispersion, skewness, kurtosis}

In this section, we show the contour plot for the index of dispersion, skewness and kurtosis 
to see the role of parameters of MCMPB$_n(\alpha,\beta,{\rm e}^{\psi})$.
The plot range is $\alpha>0$ and $\beta>0$ where the MCMPB is unimodal.
Here, we let $\mu$ be the mean of $X$ and $\mu_i$ be the $i$th central moment of $X$.

\vskip 3mm

\noindent {\bf 3.6.1. Index of dispersion}

\begin{figure}[!htb]
\begin{center}
\begin{tabular}{ccc}
{\footnotesize $n=5, \psi=-1$} & {\footnotesize $n=5, \psi=0$} & {\footnotesize $n=5, \psi=1$}
\\
\includegraphics[width=0.3\hsize,clip]{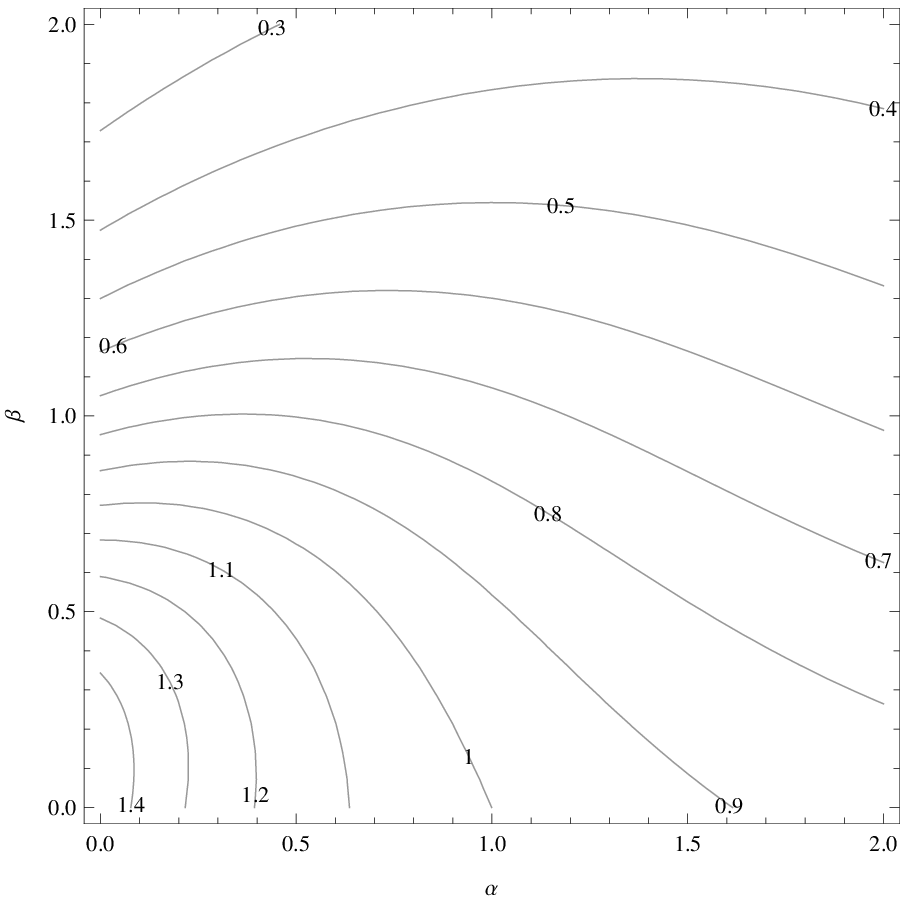}
&
\includegraphics[width=0.3\hsize,clip]{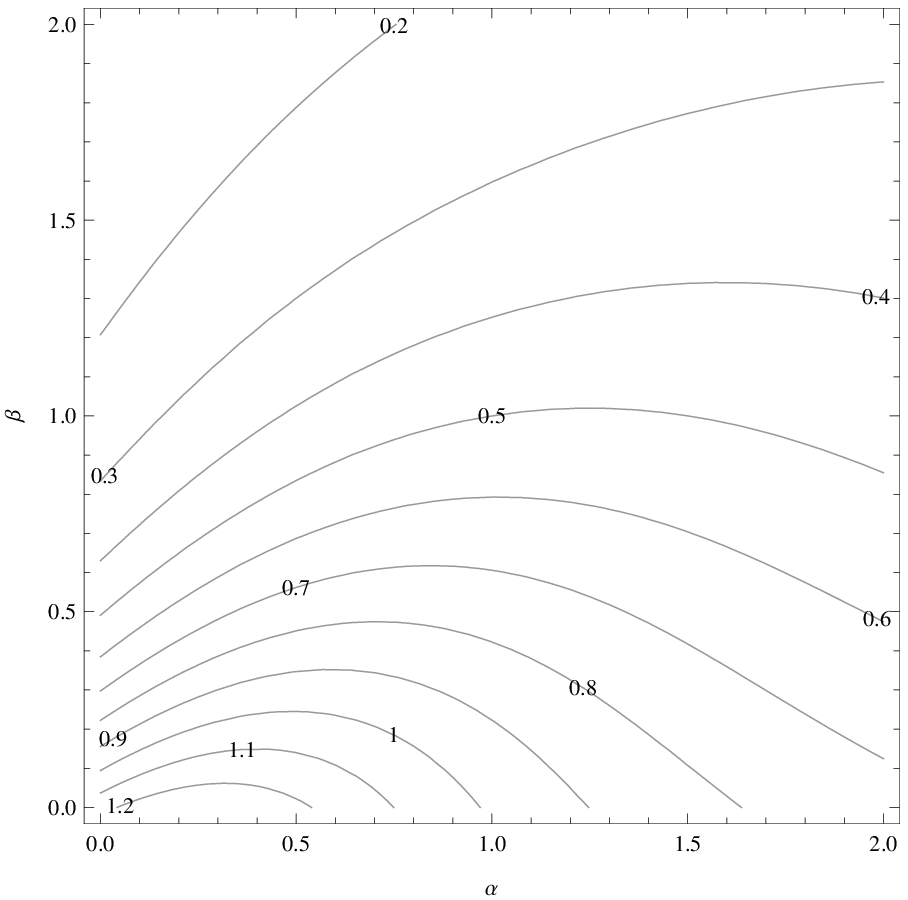}
&
\includegraphics[width=0.3\hsize,clip]{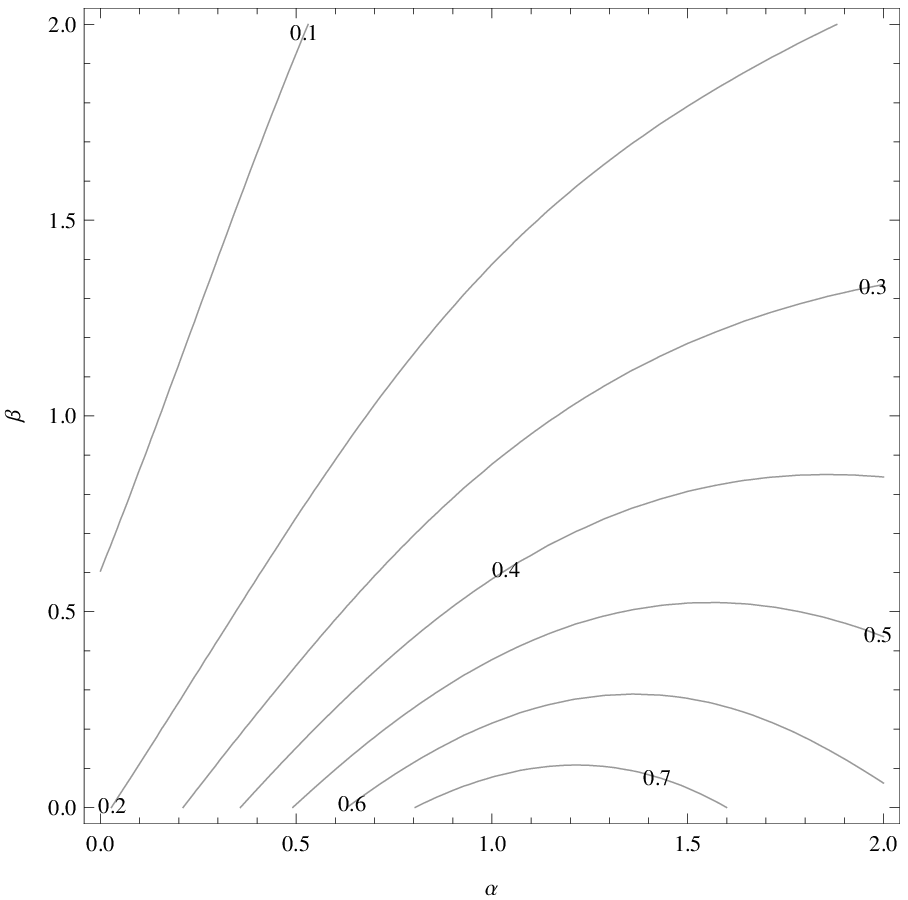}
\\
{\footnotesize $n=15, \psi=-1$} & {\footnotesize $n=15, \psi=0$} & {\footnotesize $n=15, \psi=1$}
\\
\includegraphics[width=0.3\hsize,clip]{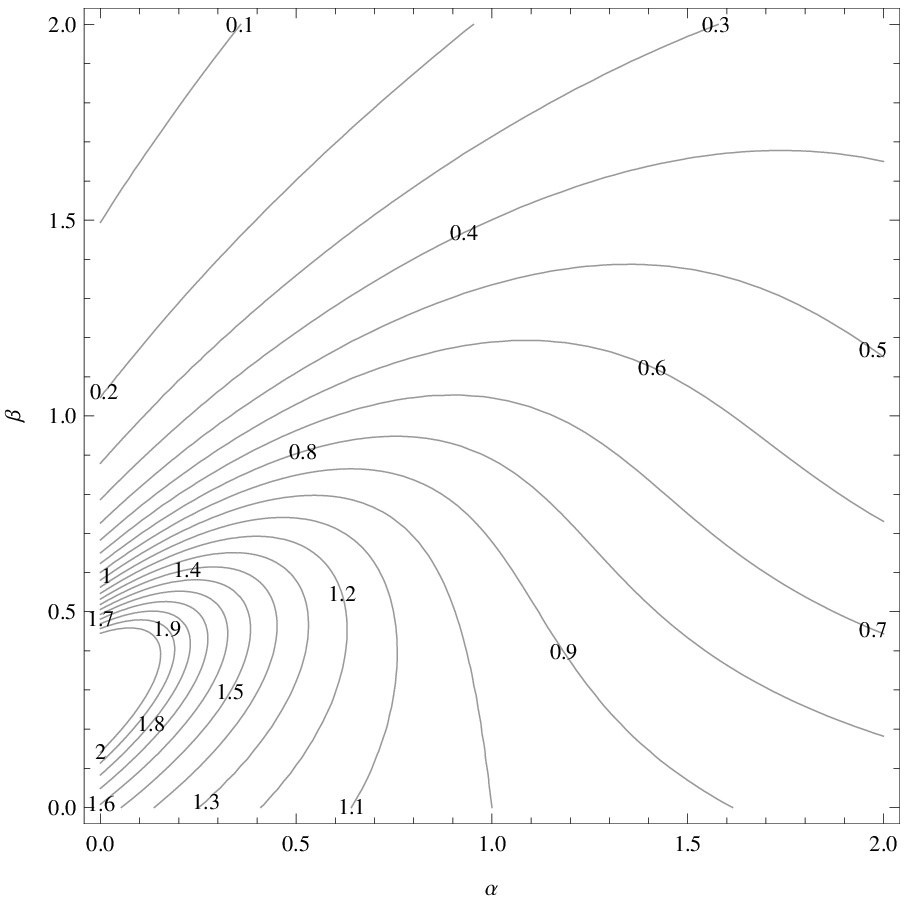}
&
\includegraphics[width=0.3\hsize,clip]{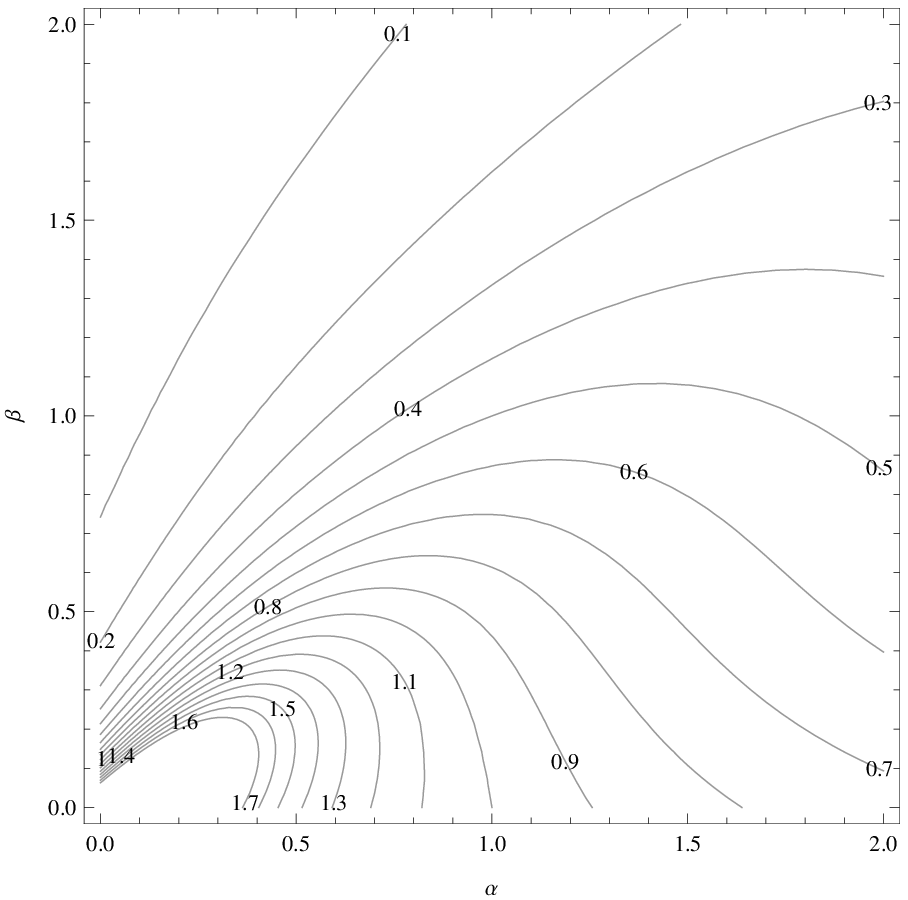}
&
\includegraphics[width=0.3\hsize,clip]{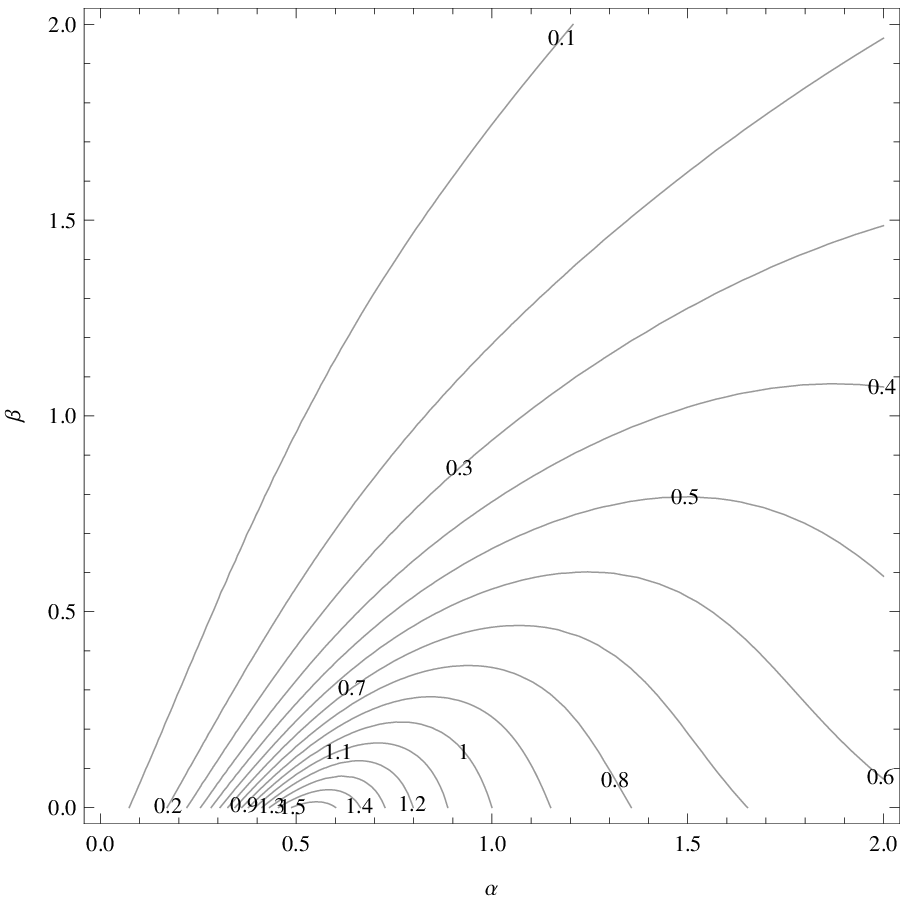}
\end{tabular}
\caption{The contour plot of dispersion.}
\label{dispersion}
\end{center}
\end{figure}

The index of dispersion is a normalized measure of the dispersion, defined by $\mu_2/\mu$.
Figure \ref{dispersion} shows the contour plot of the dispersion of the MCMPB distribution.

We can see that the region of under-dispersion becomes large as $\alpha$, $\beta$ or $\psi$ becomes large or as $n$ becomes small.
Note that the MCMPB distribution is over-dispersed for $\alpha>1$ and under-dispersion  for $\alpha<1$ when $\beta=0$.
This is because the MCMPB distribution reduces to the truncated CMP distribution when $\beta=0$.

\vskip 3mm

\noindent {\bf 3.6.2. Index of skewness}

\begin{figure}[!htb]
\begin{center}
\begin{tabular}{ccc}
{\footnotesize $n=5, \psi=-1$} & {\footnotesize $n=5, \psi=0$} & {\footnotesize $n=5, \psi=1$}
\\
\includegraphics[width=0.3\hsize,clip]{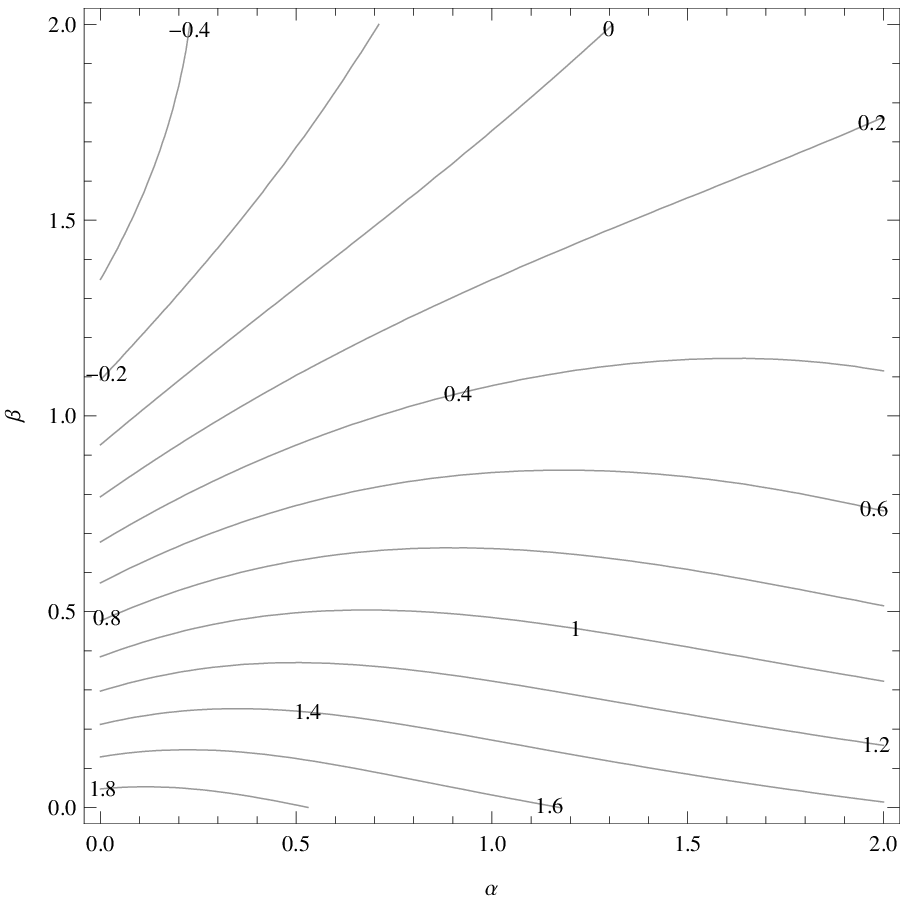}
&
\includegraphics[width=0.3\hsize,clip]{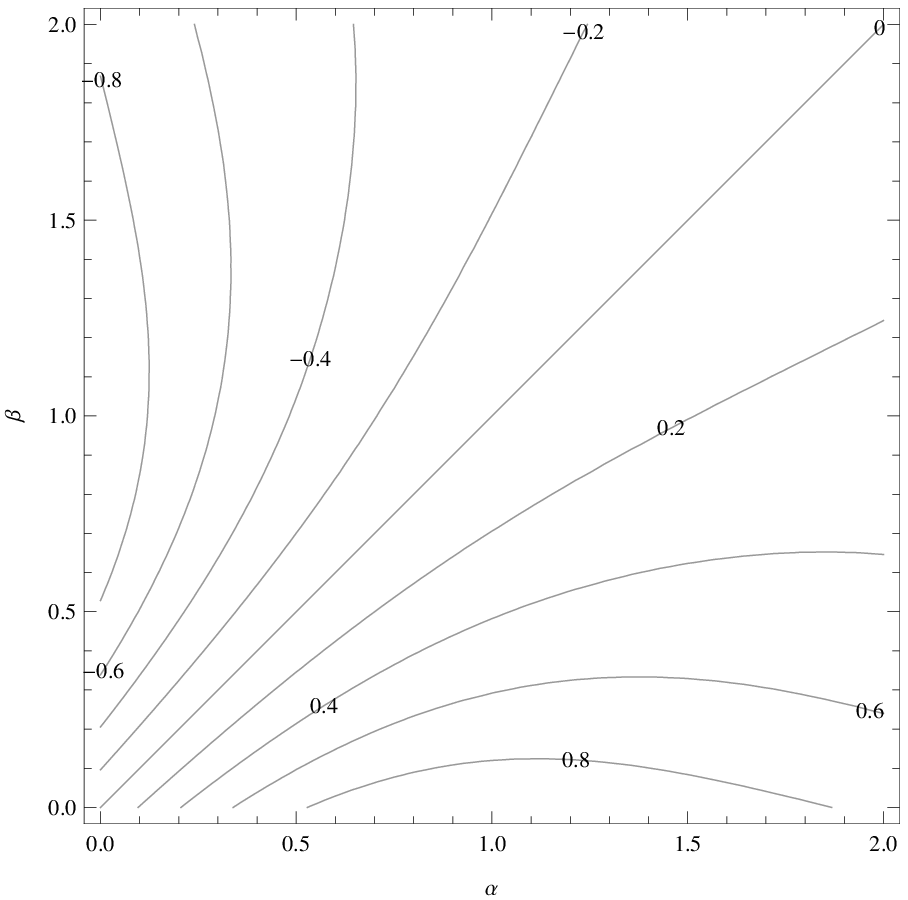}
&
\includegraphics[width=0.3\hsize,clip]{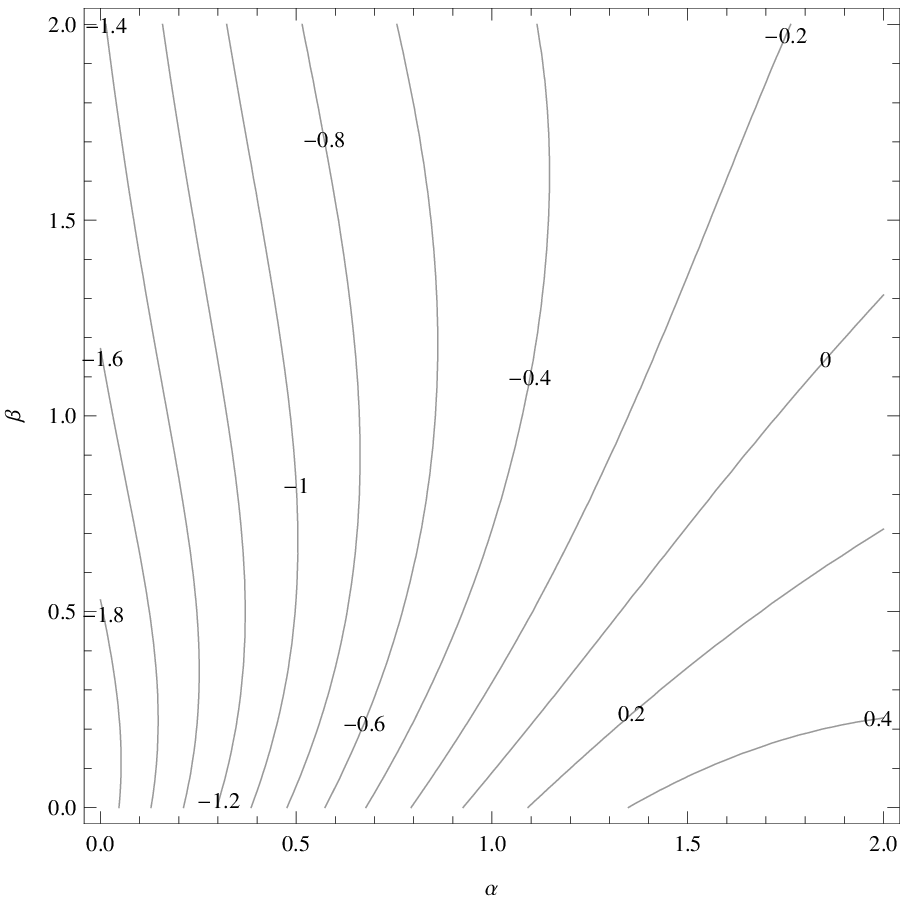}
\\
{\footnotesize $n=15, \psi=-1$} & {\footnotesize $n=15, \psi=0$} & {\footnotesize $n=15, \psi=1$}
\\
\includegraphics[width=0.3\hsize,clip]{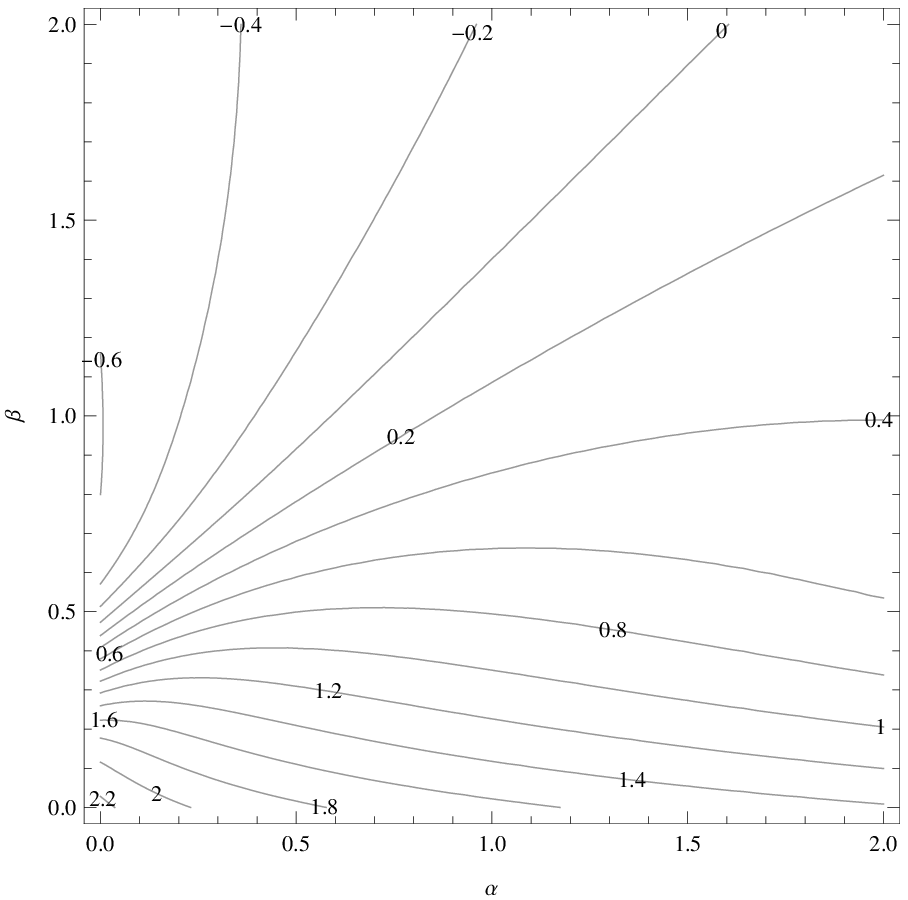}
&
\includegraphics[width=0.3\hsize,clip]{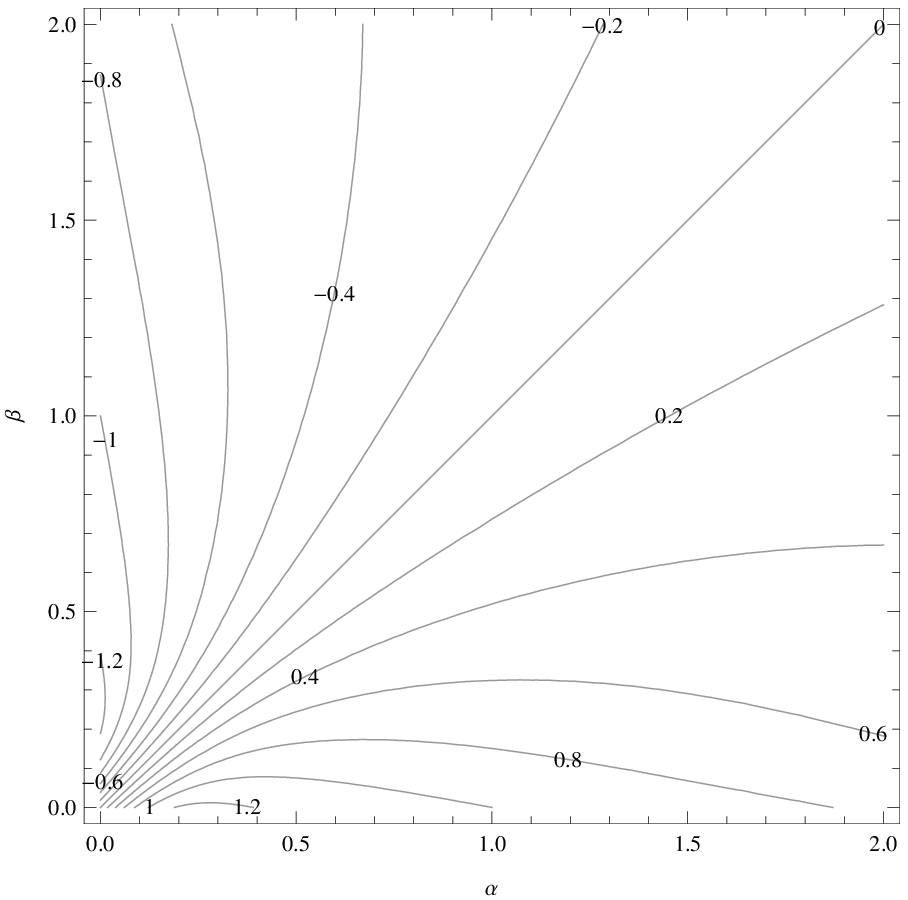}
&
\includegraphics[width=0.3\hsize,clip]{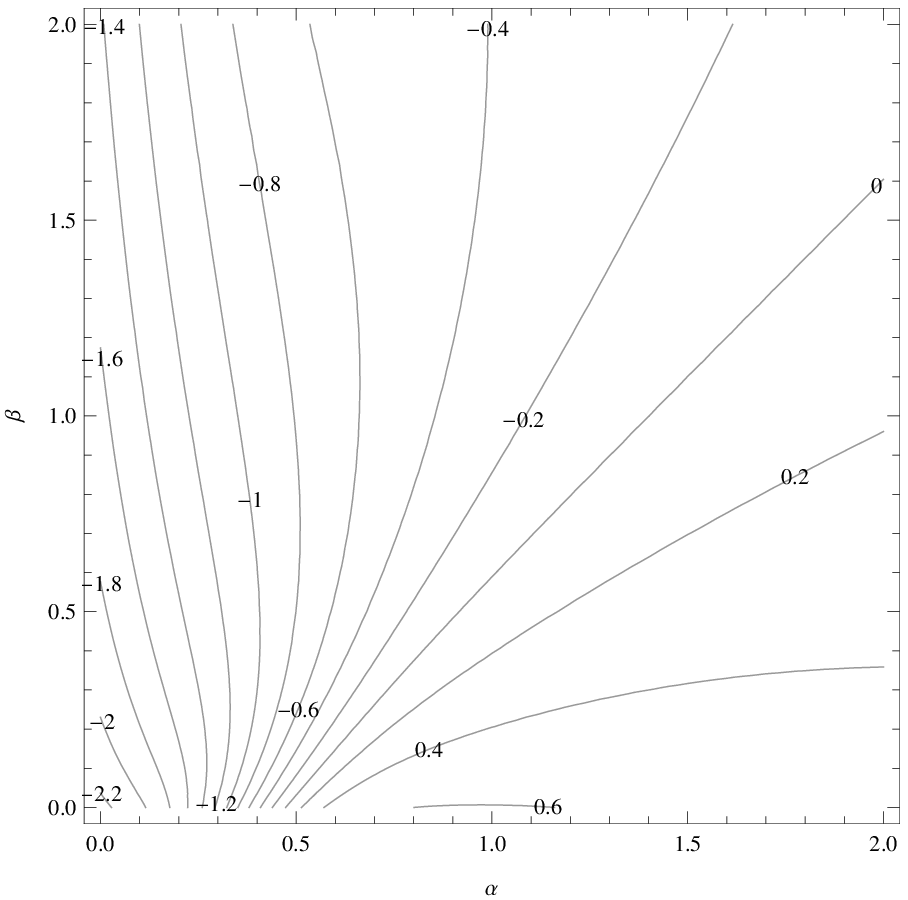}
\end{tabular}
\caption{The contour plot of skewness.}
\label{skew}
\end{center}
\end{figure}

The index of skewness is a measure of the asymmetry, defined by $\mu_3/\mu_2^{3/2}$. 
Figure \ref{skew} shows the contour plot of the skewness.

From the last statement of Section 3.4,
letting $s_n(\alpha,\beta,{\rm e}^{\psi})$ be the skewness of MCMPB$_n(\alpha,\beta,{\rm e}^{\psi})$,
we see that $s_n(\alpha,\beta,{\rm e}^{\psi})=-s_n(\beta,\alpha,{\rm e}^{-\psi})$.
This fact is confirmed from the plot of $\psi=\pm 1$.
For $\psi=0$, the MCMPB distribution is symmetry for $\alpha=\beta$, positive-skewed for $\alpha>\beta$ and negative-skewed for $\alpha<\beta$.

\vskip 3mm

\noindent {\bf 3.6.3. Index of kurtosis}

\begin{figure}[!htb]
\begin{center}
\begin{tabular}{ccc}
{\footnotesize $n=5, \psi=-1$} & {\footnotesize $n=5, \psi=0$} & {\footnotesize $n=5, \psi=1$}
\\
\includegraphics[width=0.3\hsize,clip]{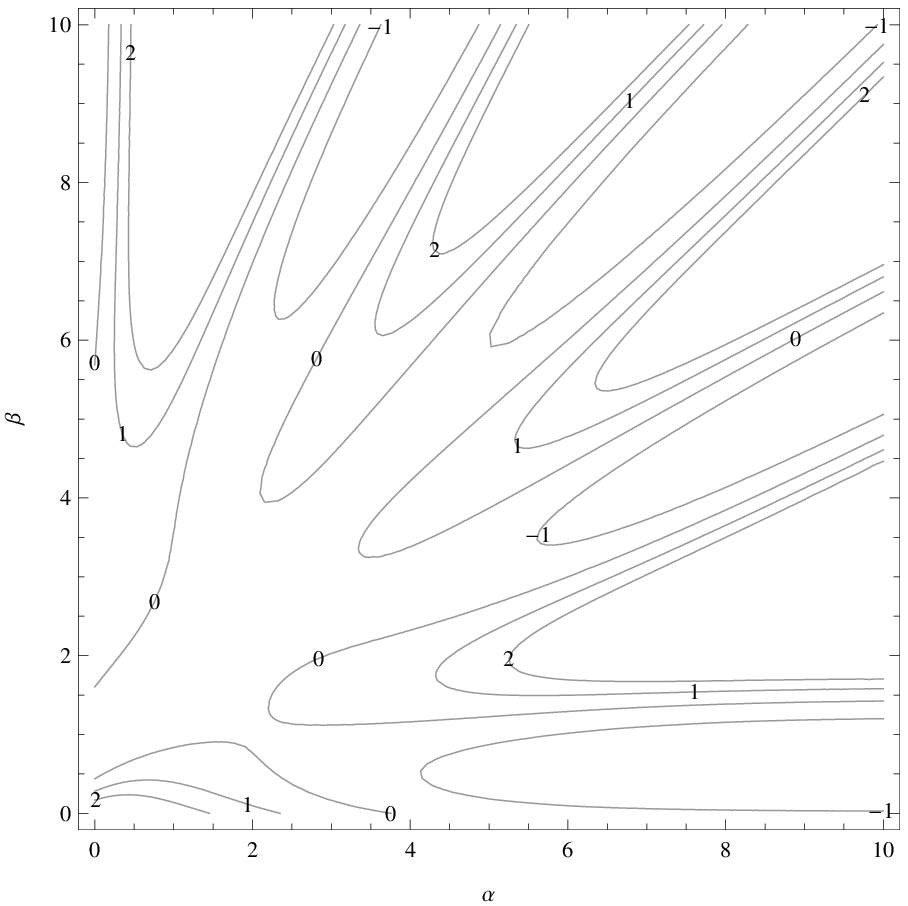}
&
\includegraphics[width=0.3\hsize,clip]{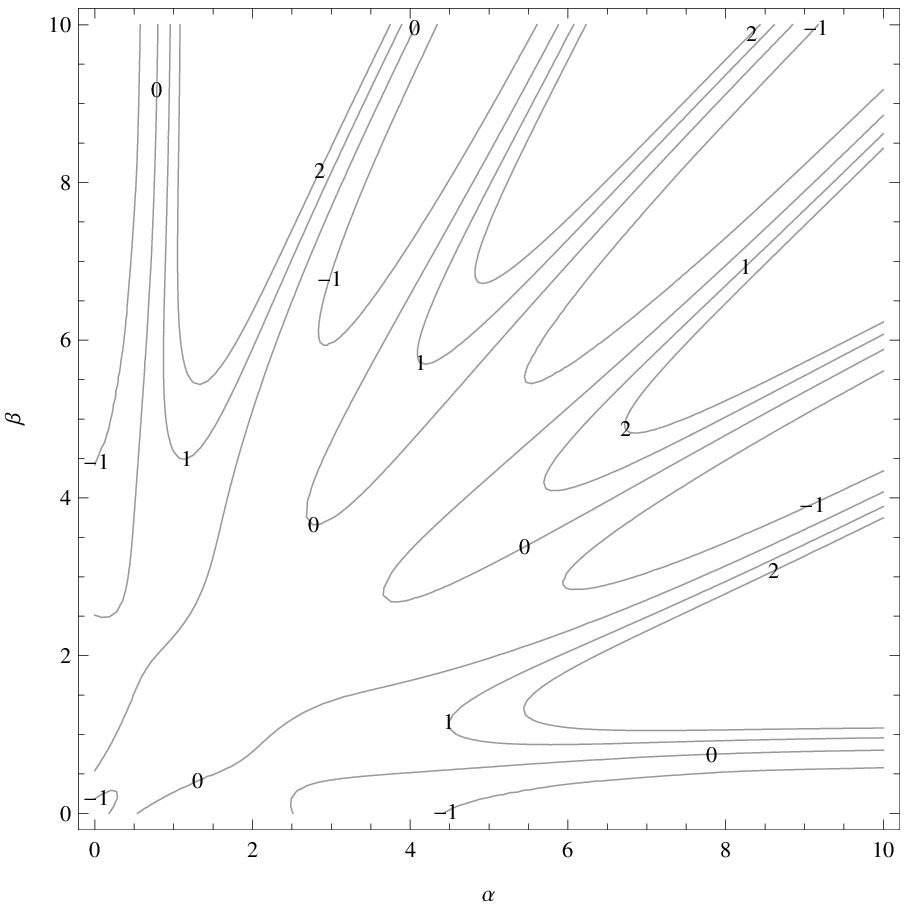}
&
\includegraphics[width=0.3\hsize,clip]{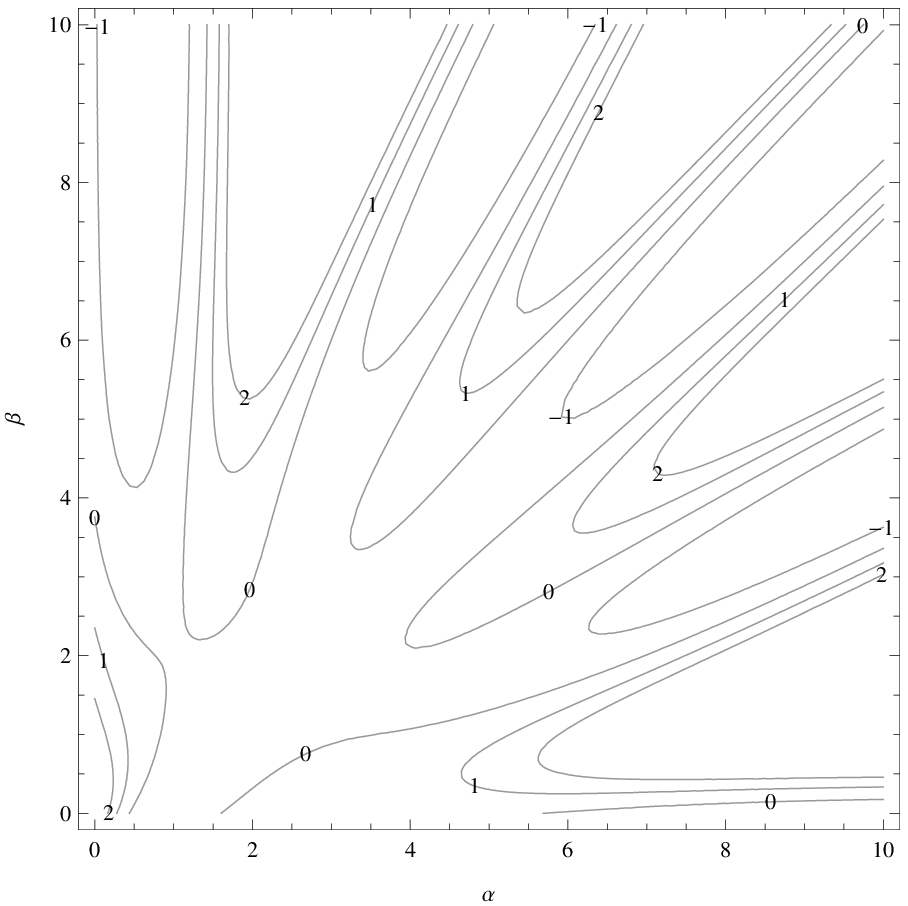}
\\
{\footnotesize $n=15, \psi=-1$} & {\footnotesize $n=15, \psi=0$} & {\footnotesize $n=15, \psi=1$}
\\
\includegraphics[width=0.3\hsize,clip]{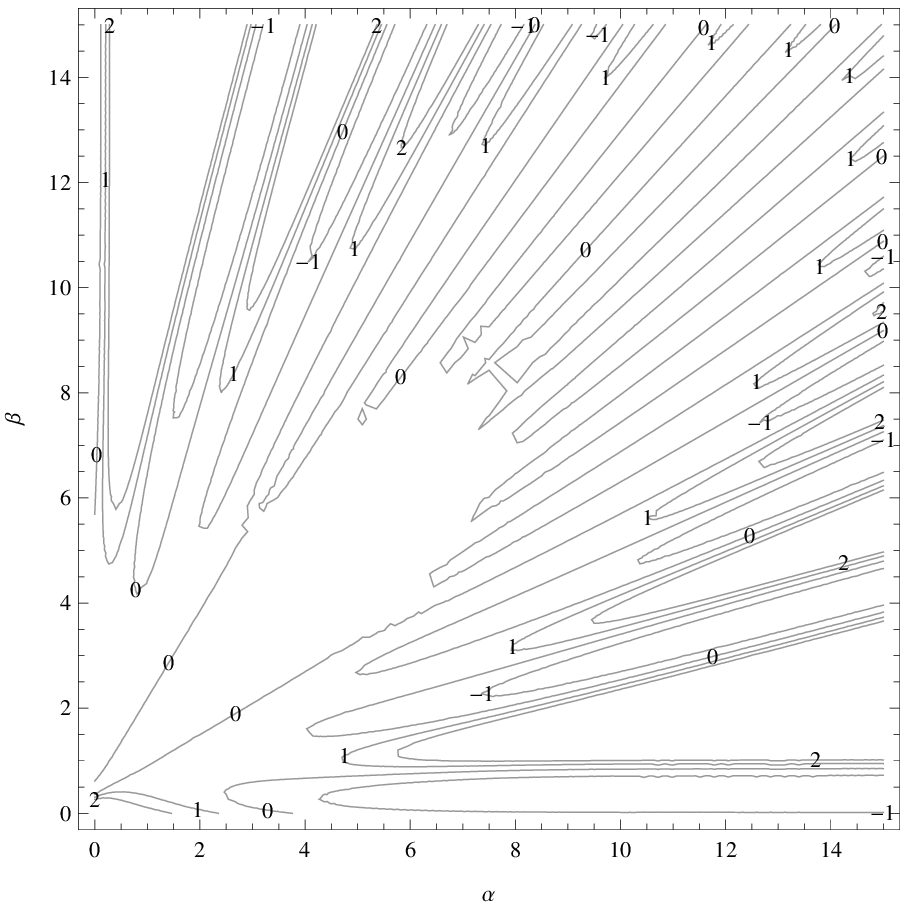}
&
\includegraphics[width=0.3\hsize,clip]{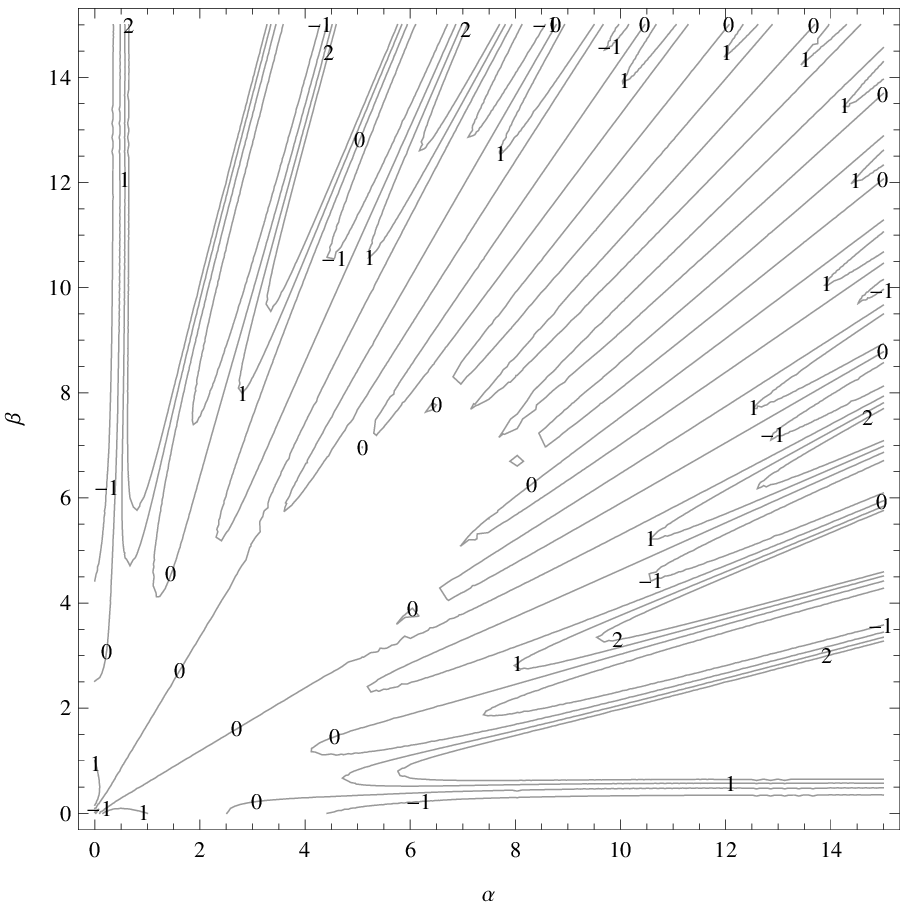}
&
\includegraphics[width=0.3\hsize,clip]{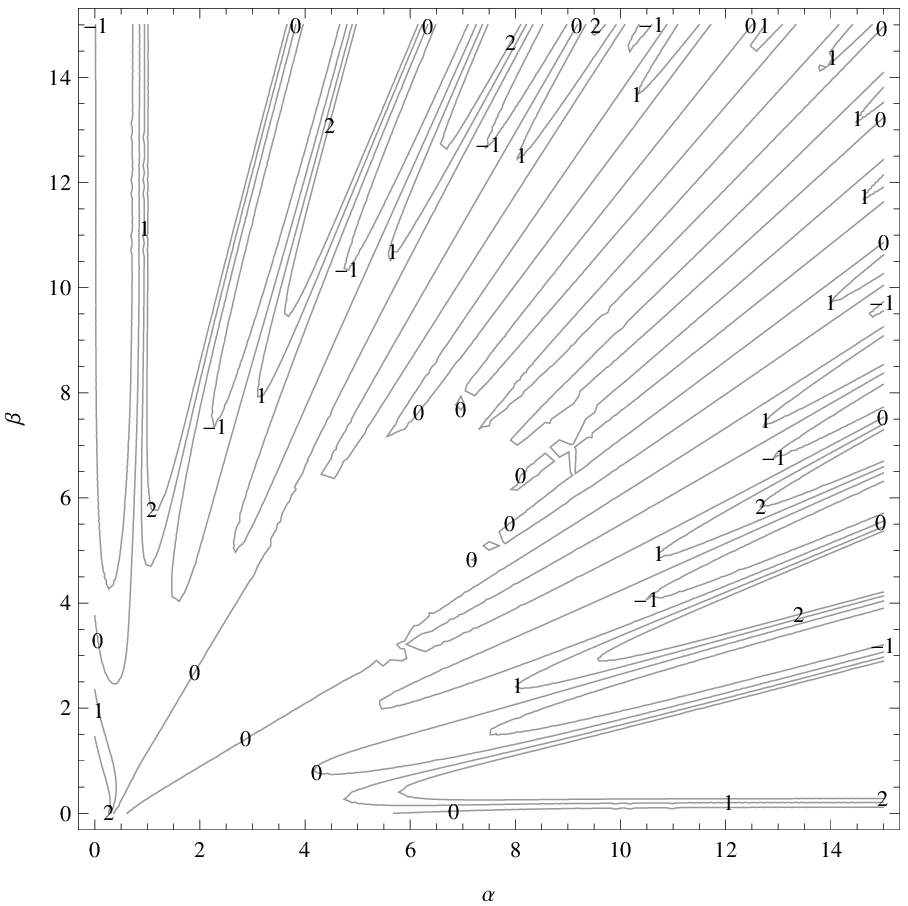}
\end{tabular}
\caption{The contour plot of kurtosis.}
\label{kurt}
\end{center}
\end{figure}

The index of kurtosis is  a measure of the tailedness, defined by $\mu_4/\mu_2^2-3$.
Figure \ref{kurt} shows the contour plot of the kurtosis of the MCMPB distribution.

From the last statement of Section 3.4,
letting $k_n(\alpha,\beta,{\rm e}^{\psi})$ be the kurtosis of MCMPB$_n(\alpha,\beta,{\rm e}^{\psi})$,
we see that $k_n(\alpha,\beta,{\rm e}^{\psi})=k_n(\beta,\alpha,{\rm e}^{-\psi})$.
This fact is confirmed from the plot of $\psi=\pm 1$.
We can see that the MCMPB distribution is platy-kurtic, or negative kurtosis, for small $\alpha$ and $\beta$
and lepto-kurtic, or positive kurtosis, for large $\alpha$ or $\beta$.
As $n$ becomes large, the range for lepto-kurtic becomes small.
Note that we can see the negative kurtosis for large $\alpha$ and $\beta$. 
This is because the MCMPB distribution can become a uniform distribution on consecutive two points
as shown in Section 3.5.

\vskip 3mm

\noindent {\bf 3.7. The other properties}

\noindent {\it Theorem 2}\\
For any bounded function $f(\cdot)$ such that the expectation exists,
$X \sim$MCMPB$_n(\alpha, \beta, \theta)$ if and only if ${\rm E}[\theta (n-X)^\beta f(X+1)-X^\alpha f(X)]=0$. \\

This theorem can be proven by the same argument by  Brown an Xia (2001).

Let us denote by $X^{(w)}$ the $w$-power biased random variable corresponding to $X$.
Then the pmf of $X^{(w)}$ is given by
\begin{eqnarray*}
{\rm P}(X^{(w)}=x)=\frac{x^w {\rm P}(X=x)}{{\rm E}[X^w]}.
\end{eqnarray*}
In particular, for $w=1$, we get size biased distribution.\\

Using Theorem 2, we can get the following theorem.

\noindent {\it Theorem 3.}\\
If $X \sim$MCMPB$_n(\alpha, \beta, \theta)$ with $\beta \geq 0$, 
then ${\rm E}[f(X^{\alpha})] \leq {\rm E}[f(X+1)]$ for all increasing function $f(\cdot)$, or
$X^{(\alpha)}$ is stochastically smaller than $X+1$.\\

\noindent {\it Proof.}
With a simple calculation, we see that
\begin{eqnarray*}
{\rm E}[f(X+1)]-{\rm E}[f(X^{(\alpha)})]= {\rm E} \left[ f(X+1)\left( 1- \frac{\theta(n-X)^{\beta}}{{\rm E}[X^{\alpha}]}\right)\right].
\end{eqnarray*}
From Theorem 2,  it is satisfied that
${\rm E} \left[1- \frac{\theta(n-X)^{\beta}}{{\rm E}[X^{\alpha}]}\right]=0$
and thus,
\begin{eqnarray*}
{\rm E} \left[ f(X+1)\left( 1- \frac{\theta(n-X)^{\beta}}{{\rm E}[X^{\alpha}]}\right)\right]
= {\rm Cov} \left[f(X+1), 1- \frac{\theta(n-X)^{\beta}}{{\rm E}[X^{\alpha}]}\right].
\end{eqnarray*}
Since $f(x+1)$ and $1- \frac{\theta(n-x)^{\beta}}{{\rm E}[X^{\alpha}]}$ are increasing functions, the above covariance is non-negative.
Therefore, we can get the result 
\begin{eqnarray*}
{\rm E}[f(X+1)] \geq {\rm E}[f(X^{(\alpha)})].
\end{eqnarray*}
\begin{flushright}$\Box$\end{flushright}

\vskip 3mm

\noindent {\bf 4. Computation of probabilities}

For the CMP distribution the computation of the normalizing constant can be involved. 
A discussion of the computational aspect is given in a recent paper of Gupta et al. (2014).

Since the pmf and the normalizing constant $C_n (\alpha,\beta,\theta)$ of the MCMPB distribution involves factorial functions 
which in numerical computations overflow very quickly, we facilitate the computation by following the approach in  Lee et al. (2001).
We derive recurrence formula to compute the normalizing constant and probabilities P$(X=x)$ of MCMPB$_n(\alpha, \beta, \theta)$.
For $C_n (\alpha,\beta,\theta)$, let
$$
C_n(\alpha,\beta,\theta)=\sum_{k=0}^n a_k,
$$
where $a_k=\frac{\theta^k}{k!^{\alpha}(n-k)!^{\beta}}$. We have
$$
\frac{a_{k+1}}{a_k}=\frac{\theta (n-k)^{\beta}}{(k+1)^{\alpha}}
$$
with $\ a_0=n!^{-\beta}$. For numerical stability, $a_0$ may be scaled to 1.

Similarly, we have
$$
\frac{{\rm P}(X=x+1)}{{\rm P}(X=x)}=\frac{\theta (n-x)^{\beta}}{(x+1)^{\alpha}}
$$
with P$(X=0)=n!^{-\beta}/C_n(\alpha,\beta,\theta))$.
The two-term recurrence formulas for $a_k$ and P$(X=x)$ avoid the computation of the factorials.

\vskip 3mm

\noindent {\bf 5. Estimation and simulation}

Let $f_i$ be the observed frequency for $i$ events and $N=\sum_{i=0}^n f_i$ be the sample size in some dataset.
Then the likelihood function of the MCMPB distribution for the dataset is given by
\begin{eqnarray}
\label{loglikelihood}
L(\alpha,\beta,\phi)=\exp[N\{\phi S_1+\alpha S_2+\beta S_3-\log C^*_n(\alpha,\beta,\phi)\}],
\end{eqnarray}
where $S_1=\sum_{i=0}^{n} i f_i/N$, $S_2=-\sum_{i=0}^{n}f_i\log i!/N$ and $S_3=-\sum_{i=0}^{n} f_i \log(n-i)!/N$.
Since the MCMPB distribution is a member of the exponential family, we see that $(S_1, S_2, S_3)$ is the minimal sufficient statistics for $(\phi, \alpha, \beta)$.
The likelihood equations are given by
\begin{eqnarray}
\label{likelihoodeq}
{\rm E}[X]=S_1,\ {\rm E}[\log X!]=S_2,\ {\rm E}[\log (n-X)!]=S_3.
\end{eqnarray}
The solution of this equations is, if exists, unique and, by allowing $(\alpha, \beta, \phi)=[-\infty, \infty]^3$, the ML estimates are uniquely determined.
Since the equations (\ref{likelihoodeq}) cannot be solved analytically, 
an iterative method such as the Nelder--Mead method for maximizing the function (\ref{loglikelihood}) or 
the Newton--Raphson or Score method for solving the equation (\ref{likelihoodeq}) is necessary.
In this paper, we maximize the function (\ref{loglikelihood}) by using {\it NMaximize} command with the Nelder--Mead method in Mathematica 8.0.

The Fisher information matrix $I(\phi,\alpha,\beta)$ is given by {\bf Var}$^{-1}(X, \log x!, \log (n-X)!)$, 
the variance-covariance matrix of the random vector $(X, \log x!, \log (n-X)!)$.
This can be calculated by the equations (\ref{covariance}).
Since the ML estimates $(\hat{\phi},\hat{\alpha},\hat{\beta})$ are uniquely determined,
these are asymptotically distributed as the trivariate normal distribution with the mean $(\phi,\alpha,\beta)$ and the variance $I^{-1}(\phi,\alpha,\beta)$.
From this fact, we can easily construct the $95\%$ confidence intervals for the parameters $(\phi,\alpha,\beta)$.

A simulation study is conducted with $10,000$ Monte Carlo repetitions for MCMPB$_{15}(0.2, 0.4, {\ e}^{0})$ (under-dispersed case), 
MCMPB$_{15}(0.5,0.2,{\rm e}^{0.5})$ (over-dispersed case) and MCMPB$_{15}(-0.5, 0.7,{\rm e}^{-2.4})$ (bimodal case). 
The sample size $N$ is set at $100$, $500$ and $1,000$ to represent small, medium and large sample sizes.
The results are in Table \ref{simulation} with the bias, mean squared error (MSE) and the number of the trials where the true parameters are in $95\%$ confidence intervals.

\begin{table}[htbp]
\caption{Simulations for the ML estimates.}
\label{simulation}
\begin{center}
{\scriptsize \begin{tabular}[c]{ll|ccc|ccc|ccc}
     & $\alpha=0.2$ & $\beta=0.4$ & $\phi=0$  & $\alpha=0.5$ & $\beta=0.2$ & $\phi=0.5$  & $\alpha=-0.5$ & $\beta=0.7$ & $\phi=-2.4$   \\
\hline
N=100 &Bias              & $0.19$ & $-0.06$ & $0.55$ & $0.00$ & $0.05$ & $-0.11$ & $-0.01$ & $0.03$ & $-0.07$ \\
            &MSE              & $0.41$ & $0.08$ & $3.67$ & $0.09$ & $0.12$ & $1.60$ & $0.03$ & $0.05$ & $0.61$ \\
            &Number     & $9557$  & $9540$ & $9552$ & $9540$ & $9574$ & $9571$ & $9522$ & $9511$ & $9516$ \\
\hline
N=500 &Bias              & $0.03$ & $-0.01$ & $0.10$ & $0.00$ & $0.01$ & $-0.01$ & $0.00$ & $0.00$ & $-0.01$ \\
            &MSE              & $0.06$ & $0.01$ & $0.53$ & $0.02$ & $0.02$ & $0.29$ & $0.01$ & $0.01$ & $0.11$ \\
            &Number    & $9529$  & $9542$ & $9534$ & $9463$ & $9493$ & $9463$ & $9492$ & $9501$ & $9506$ \\
\hline
N=1000 &Bias              & $0.02$ & $-0.01$ & $0.05$ & $0.00$ & $0.01$ & $-0.01$ & $0.00$ & $0.00$ & $-0.01$ \\
              &MSE              & $0.03$ & $0.01$ & $0.26$ & $0.01$ & $0.01$ & $0.14$ & $0.00$ & $0.00$ & $0.06$ \\
              &Number  & $9529$  & $9519$ & $9525$ & $9503$ & $9526$ & $9518$ & $9499$ & $9503$ & $9513$ \\
\hline
\end{tabular}
}
\end{center}
\end{table}

\vskip 3mm

\noindent {\bf 6. Illustrative examples}

The application of the proposed MCMPB distribution is further illustrated using two dataset. 
These fitting results are compared with CMPB (MCMPB$_n(\alpha,\alpha,\phi)$), Beta-binomial (BB) 
$$
{\rm P}(X=x)={n \choose x} \frac{B(a+x,n+b-x)}{B(a, n+b)},
$$
negative binomial (NB) 
$$
{\rm P}(X=x)={r+x-1 \choose x} p^x(1-p)^r
$$
and CMP (\ref{compmf}) distributions. 
The parameters of these distributions are estimated by the maximum likelihood estimation.
For MCMPB, CMPB and BB distributions, we use the profile maximum likelihood method for determining the parameter $n$. 

The First dataset is sharp-topped data on the number of eggs per nest of linnets (Heyde and Schuh 1978).
In this data, the observations are counted from $1$, so we use the zero-truncated version of the distributions.
As we can see, the fitting result by MCMPB distribution is improved well in the sense of AIC and $\chi^2$ statistics.

\begin{table}[htbp]
\caption{Observed and expected frequencies of clutch size data for linnets.}
\label{linnet}
\begin{center}
{\scriptsize \begin{tabular}[c]{ccccccc}
Count & Observed & MCMPB & CMPB & BB & CMP & NB\\
\hline
1 & 18   & 24.26     & 0.00      & 32.50    & 0.00      & 243.51\\
2 & 35   & 25.28     & 2.39      & 198.45  & 0.75      & 565.46\\
3 & 210 & 175.13   & 170.74      & 673.24  & 160.16  & 875.70\\
4 & 1355& 1458.22& 1820.16  & 1370.38& 1970.37& 1017.48\\
5 & 3492& 3392.79& 2898.81 1673.65& 2662.10& 946.13\\
6 & 299  & 338.29  & 517.93& 1135.57& 591.66  & 733.41\\
7 & 5      & 0.03      & 3.95  & 330.21  & 28.59    & 487.48\\
\hline
AIC             & & 10615.16 & 11299.48 & 14179.45 &11834.34 & 18971.26 \\
$\chi^2$   & & 15.99 & 284.87 & 3372.15 & 679.90 & 8909.85 \\
d.f.             & & 2 & 1 & 4 & 4 & 2 \\
$p$-value  & & 0.00 & 0.00 & 0.00 & 0.00 & 0.00 \\
\hdashline
Parameters& Mean$=4.70$        &  $\hat{\alpha}= -10.24$ &$\hat{\alpha}= 3.72$ & $\hat{a} = 1.43 \times 10^7 $  & $\hat{r} = 9.90$              & $\hat{r} = 2750.97$ \\
& Variance$=0.48$   &  $\hat{\beta}= 12.37$    &$\hat{\psi}= 2.37$     & $\hat{b} = 7.05 \times 10^6$ & $\hat{\lambda} = 1.12 \times 10^7$ & $\hat{p} = 0.001$ \\
& Skewness$=-1.24$ &  $\hat{\psi}= -29.22$      &                                     &                                 &                                          &                                 \\
& Kurtosis$=3.56$ & & & & & \\
\hline
\end{tabular}}
\end{center}
\end{table}

The second dataset is over-dispersed data on the number of trips made by Dutch households
owning at least one car during a particular survey week in 1989 (van Ophem 2000). 
In van Ophem (2000), all proposed models to fit this slightly over-dispersed data set 
have been strongly rejected at the $0.01$ level of significance. 
Clearly from the $\chi^2$ statistics as well as the $p$-values, the MCMPB distribution provides a better
fit to this data set and will not be rejected at the 0.01 level of significance.

 \begin{table}[htbp]
\caption{Observed and expected frequencies of trips made by Dutch households owning at least one car during a particular survey week in 1989.}
\label{trip}
\begin{center}{\scriptsize
\begin{tabular}[c]{ccccccc}
Count & Observed & MCMPB & CMPB & BB & CMP & NB\\
\hline
0 & 75   &  81.24 & 112.77 & 105.93 & 100.05 & 102.43 \\
1 & 312 & 282.12& 271.93 & 278.95 & 275.11 & 281.46 \\
2 & 384 & 426.36& 384.62 & 392.98 & 398.96 & 400.11 \\
3 & 421 & 410.80& 390.08 & 389.23 & 397.93 & 391.91 \\
4 & 307 & 296.12& 307.35 & 300.96 & 304.35 & 297.26 \\
5 & 183 & 174.89& 196.22 & 191.43 & 189.44 & 186.04 \\
6 & 77   & 89.96  & 104.04 & 103.10 & 99.65   & 99.99   \\
7 & 47   & 42.18  & 46.51   & 47.79   & 45.47   & 47.43   \\
8 & 15   & 18.72  & 17.68   & 19.21   & 18.34   & 20.25   \\
9 & 9     & 8.12    & 5.74     & 6.71     & 6.64     & 7.90     \\
10& 5    & 3.56    & 1.59     & 2.03     & 2.18     & 2.85     \\
11& 0    & 1.63    & 0.38     & 0.53     & 0.65     & 0.96     \\
12& 0    & 0.81    & 0.07     & 0.12     & 0.18     & 0.30     \\
13& 1    & 0.45    & 0.01     & 0.02     & 0.05     & 0.09     \\
14& 2    & 0.31    & 0.00     & 0.00     & 0.01     & 0.03     \\
15& 0    & 0.27    & 0.17     & 0.00     & 0.00     & 0.01     \\
16& 0    & 0.37    & 0.00     & 0.24     & 0.00     & 0.00     \\
17& 1    & 1.95    & 0.00     & 0.00     & 0.00     & 0.00    \\
\hline
AIC             & & 7194.30 & 7265.32 & 7252.61 & 7233.06 &  7224.20 \\
$\chi^2$   & & 12.05 & 41.97 & 30.97 & 26.24 & 23.32 \\
d.f.             & & 6 & 7 & 7 & 7 & 7 \\
$p$-value  & & 0.06 & 0.00 & 0.00 & 0.00 & 0.00 \\
\hdashline
Parameters& Mean$=3.04$        &  $\hat{\alpha}= 1.31$ &$\hat{\alpha}= 0.71$ & $\hat{a} = 8.59$  & $\hat{r} = 0.92$              & $\hat{r} = 28.80$ \\
& Variance$=3.41$   &  $\hat{\beta}= -1.26$ &$\hat{\psi}= -1.13$   & $\hat{b} = 39.48$ & $\hat{\lambda} = 2.75$ & $\hat{p} = 0.095$ \\
& Skewness$=1.14$ &  $\hat{\psi}= 4.81$      &                                     &                                 &                                          &                                 \\
& Kurtosis$=3.62$ & & & & & \\
\hline
\end{tabular}}
\end{center}
\end{table}

\vskip 3mm

\noindent {\bf 7. Concluding remarks}

In this paper, we proposed the modified CMPB distribution 
which may arise from the conditional CMP distribution given a sum of two CMP variables or 
from a finite capacity queueing system with state-dependent arrival and service rates.
The advantage of this distribution is its flexibility of the dispersion, skewness and kurtosis and modality.
The normalizing constant is not a closed form, but it includes only the summation of finite series and thus, 
need not any special approximation.
Since the proposed distribution belongs to the exponential family, we can easily drive properties such as the moments and ML estimates.
From the geneses of the MCMPB distribution, fitting the MCMPB distribution to real dataset gives some information about the dataset as described 
in examples of Section 2.
In this paper, we used only ML estimation, but the Baysian inference will be also available for fitting through the conjugate prior distribution.
It is difficult to find the roles of parameters, or which parameters control the mean, variance and skewness of the MCMPB distribution.
However, it will be easy to find the roles for the special cases of the MCMPB distribution. 
In this sense, it is valuable to study and discuss about the special cases of the MCMPB distribution.

\vskip 3mm

\noindent {\bf Acknowledgments}

This research is partially supported by ISM FY2014 travel support for an international paper/poster presentation by young researchers. 
(ISM stands for The Institute of Statistical Mathematics in Japan)

\vskip 3mm

\noindent {\bf Bibliography}
\vskip 3mm

\noindent Altham, P. M. E. (1978). 
Two generalizations of the binomial distribution. 
{\it Journal of the Royal Statistical Society, Series C (Applied Statistics)}, {\bf 27}, 162--167.

\vskip 3mm

\noindent Bliss, C. I. (1953).
Fitting the negative binomial distribution to biological data.
{\it Biometrics}, {\bf 9}, 176--200.

\noindent Bahadur R. R. (1961).
A representation of the joint distribution of responses to n dichotomous items. 
{\it In Studies in Item Analysis and Prediction (ed. H. Solomon), 
Stanford Mathematical Studies in the Social Sciences VI}, California: Stanford University Press, Stanford, 158--168.

\noindent Borges, P., Rodrigues, J., Balakrishnan, and N., Bazan, J. (2014). 
A COM-Poisson type generalization of the binomial distribution and its properties and applications. 
{\it Statistics \& Probability Letters}, {\bf 87}, 158--166.

\noindent Brown, T. C., and Xia A (2001). 
Steinfs Method and Birth-Death Processes. 
{\it The Annals of Probability}, {\bf 29}, 1373--1403.

\noindent Castillo, J., and P\'{e}rez-Casany, M. (2005).
Overdispersed and underdispersed Poisson generalizations.
{\it Journal of Statistical Planning and Inference}, {\bf 134}, 486--500.

\noindent Chang, Y., and Zelterman, D. (2002). 
Sums of dependent Bernoulli random variables and disease clustering. 
{\it Statistics \& Probability Letters}, {\bf 57}, 363--373.

\noindent Consul, P. C. (1974). 
A simple urn model dependent upon predetermined strategy.
{\it Sankhy\={a}, Series B}, {\bf 36}, 391--399.

\noindent Consul, P. C., and Famoye, F. (2006). 
{\it Lagrange Probability Distributions.}
New York: Birkh\"{a}user.

\noindent Consul, P. C., and Mittal, S. P. (1975).
A new urn model with predetermined strategy.
{\it Biometrische Zeitschrift}, {\bf 17}, 67--75.

\noindent Conway, R. W., and Maxwell, W. L. (1962). 
A queueing model with state dependent service rates.
{\it Journal of Industrial Engineering}, {\bf 12}, 132--136.

\noindent Gupta, R. C., Sim, S. Z., and Ong, S. H. (2014). 
Analysis of Discrete Data by Conway--Maxwell Poisson distribution. 
{\it AStA Advances in Statistical Analysis}, {\bf 98}, 327--343.

\noindent Heyde, C. C., and Schuh, H-J. (1978). 
Uniform bounding of probability generating functions and the evolution of reproduction rates in birds. 
{\it Journal of Applied Probability}, {\bf 15}, 243--250.

\noindent Jinkinson, R. A., and Slater, M. (1981). 
Critical discussion of a graphical method for identifying discrete distributions. 
{\it Journal of the Royal Statistical Society Series D (The Statistician)}, {\bf 30}, 239--248.

\noindent Kitano, M., Shimizu, K., and Ong, S. H. (2005).
The generalized Charlier series distribution as a distribution with two-step recursion.
{\it Statistics \& Probability Letters}, {\bf 75}, 280--290.

\noindent Lee, P. A., Ong, S. H., and Srivastava, H. M. (2001). 
Some integrals of the products of Laguerre polynomials.
{\it International Journal of Computer Mathematics}, {\bf 78}, 303--321.

\noindent Lindsey, J. K. (1995).
{\it Modelling frequency and count data.}
New York: Oxford University Press.

\noindent Ong, S. H. (1988).
A discrete Charlier series distribution.
{\it Biometrical Journal}, {\bf 30}, 1003--1009.

\noindent Rao, C. R. (1965).
On discrete distributions arising out of methods of ascertainment.
{\it Sankhy\={a}, Series A}, {\bf  27}, 311--324.

\noindent Upton, G. J. G., and Lampitt, G. A. (1981). 
A model for interyear change in the size of bird populations. 
{\it Biometrics}, {\bf 37}, 113--127.

\noindent Shmueli, G., Minka, T. P., Kadane, J. B., Borle, S., and Boatwright, P. (2005).
A useful distribution for fitting discrete data: revival of the Conway--Maxwell--Poisson distribution.
{\it Journal of the Royal Statistical Society, Series C (Applied Statistics)}, {\bf 54}, 127--142.

\noindent Sokal, R. R., and Rohlf, F. J. (1994).
{\it Biometry: The Principles and Practices of Statistics in Biological Research.}
New York: W.H. Freeman.

\noindent van Ophem, H. (2000). 
Modeling selectivity in count-data models.
{\it Journal of Business \& Economic Statistics}, {\bf 18}, 503--511.

\end{document}